\documentclass[12pt,english]{article}
\usepackage[T1]{fontenc}
\usepackage[latin9]{inputenc}
\usepackage[a4paper]{geometry}
\usepackage{babel}
\usepackage{mathtools}
\usepackage{enumitem}
\usepackage{amsmath}
\usepackage{amsthm}
\usepackage{amssymb}
\usepackage[numbers]{natbib}
\usepackage[unicode=true,pdfusetitle,
 bookmarks=true,bookmarksnumbered=false,bookmarksopen=false,
 breaklinks=false,pdfborder={0 0 1},backref=false,colorlinks=false]
 {hyperref}

\makeatletter
\theoremstyle{plain}
\ifx\thechapter\undefined
	\newtheorem{thm}{\protect\theoremname}
\else
	\newtheorem{thm}{\protect\theoremname}[chapter]
\fi
\theoremstyle{plain}
\newtheorem{prop}[thm]{\protect\propositionname}
\theoremstyle{remark}
\newtheorem{rem}[thm]{\protect\remarkname}
\theoremstyle{plain}
\newtheorem{cor}[thm]{\protect\corollaryname}
\theoremstyle{plain}
\newtheorem{lem}[thm]{\protect\lemmaname}

\usepackage{babel}
\usepackage{titlesec}
\usepackage[titletoc,toc,title]{appendix}
\usepackage{xcolor}
\usepackage{graphicx}

\theoremstyle{plain}
\newtheorem{assumption}[thm]{Assumption}

\DeclareMathOperator{\E}{{\mathbb E}}

\DeclareMathOperator{\D}{{\mathbb D}}
\DeclareMathOperator{\R}{{\mathbb R}}

\DeclareMathOperator{\N}{{\mathbb N}}

\DeclareMathOperator{\Q}{{\mathbb Q}}

\DeclareMathOperator{\dom}{dom}

\newcommand{\vertiii}[1]{{\left\vert\kern-0.25ex\left\vert\kern-0.25ex\left\vert #1 
    \right\vert\kern-0.25ex\right\vert\kern-0.25ex\right\vert}}

\makeatother

\providecommand{\corollaryname}{Corollary}
\providecommand{\lemmaname}{Lemma}
\providecommand{\propositionname}{Proposition}
\providecommand{\remarkname}{Remark}
\providecommand{\theoremname}{Theorem}

\begin{document}
\global\long\def\epsilon{\varepsilon}

\global\long\def\E{\mathbb{E}}

\global\long\def\I{\mathbf{1}}

\global\long\def\N{\mathbb{N}}

\global\long\def\R{\mathbb{R}}

\global\long\def\C{\mathbb{C}}

\global\long\def\Q{\mathbb{Q}}

\global\long\def\P{\mathbb{P}}

\global\long\def\D{\Delta_{n}}

\global\long\def\dom{\operatorname{dom}}

\global\long\def\b#1{\mathbb{#1}}

\global\long\def\c#1{\mathcal{#1}}

\global\long\def\s#1{{\scriptstyle #1}}

\global\long\def\u#1#2{\underset{#2}{\underbrace{#1}}}

\global\long\def\r#1{\xrightarrow{#1}}

\global\long\def\mr#1{\mathrel{\raisebox{-2pt}{\ensuremath{\xrightarrow{#1}}}}}

\global\long\def\t#1{\left.#1\right|}

\global\long\def\l#1{\left.#1\right|}

\global\long\def\f#1{\lfloor#1\rfloor}

\global\long\def\sc#1#2{\langle#1,#2\rangle}

\global\long\def\abs#1{\lvert#1\rvert}

\global\long\def\bnorm#1{\Bigl\lVert#1\Bigr\rVert}

\global\long\def\wraum{(\Omega,\c F,\P)}

\global\long\def\fwraum{(\Omega,\c F,\P,(\c F_{t}))}

\global\long\def\norm#1{\lVert#1\rVert}

\author{Randolf Altmeyer \\ 
\\
\textit{ Humboldt-Universit\"at zu Berlin \footnote{Institut f\"ur Mathematik, Humboldt-Universit\"at zu Berlin, Unter den Linden 6, 10099 Berlin, Germany. Email: altmeyrx@math.hu-berlin.de}}} 

\title{Approximation of occupation time functionals}

\maketitle
\begin{abstract}
The strong $L^2$-approximation of occupation time functionals is studied with respect to  discrete observations of a $d$-dimensional càdlàg process. Upper bounds on the error are obtained under weak assumptions, generalizing  previous results in the literature considerably. The approach relies on regularity for the marginals of the process and applies also to non-Markovian processes, such as fractional Brownian motion. The results are used to approximate occupation times and local times. For Brownian motion, the upper bounds are shown to be sharp up to a log-factor.

\medskip


 \end{abstract} 

\section{Introduction}

The approximation of integral-type functionals for random integrands
is a classical problem. It appears in the study of numerical approximation
schemes for stochastic differential equations (\citep{Gobet2008,Kohatsu-Higa2014,Neuenkirch:2019uo})
and in the analysis of statistical methods for stochastic processes
(\citep{Chorowski2015a,Florens-Zmirou1993,Jacod1998}). Early works
focused on choosing optimal sampling times (\citep{sacks1966designs})
or on random integrands as a tool for Bayesian numerical analysis
(cf. \citep{diaconis1988bayesian} or \citep{ritter2007average} for
an overview). Recently, there has been growing interest in estimating
integral functionals of the form
\[
\Gamma_{T}\left(f\right)=\int_{0}^{T}f\left(X_{t}\right)dt
\]
for a known measurable function $f$ and an $\R^{d}$-valued stochastic
process $X=(X_{t})_{0\leq t\leq T}$, $T>0$. Such functionals are
called \emph{occupation time functionals,} as they generalize the
occupation time $\Gamma_{T}(\I_{A})$ of a set $A\subset\R^{d}$.

Suppose we have access to $X_{t_{k}}$ at discrete time points $t_{k}=k\Delta_{n}$,
where $\Delta_{n}=T/n$ and $k\in\{0,\dots,n\}$. The paths of $f(X)$
are typically rough, even for smooth $f$, allowing only for lower
order quadrature rules to approximate $\Gamma_{T}(f)$ (\citep{Cruz-Uribe2002}).
A natural estimation scheme is given by the Riemann estimator
\[
\widehat{\Gamma}_{T,n}(f)=\Delta_{n}\sum_{k=1}^{n}f(X_{t_{k-1}}).
\]
Its theoretical properties have been considered systematically only
in few works and only for rather specific processes $X$ and functions
$f$. The goal of this paper is to study in a general setting the
strong $L^{2}(\P)$-approximation of $\Gamma_{T}(f)$ by $\widehat{\Gamma}_{T,n}(f)$
and to derive upper bounds on the error, which are explicit in terms
of $f$, $T$ and $\Delta_{n}$, and depend on the dimension $d$
only through $f$. The function $f$ is considered in Hölder or fractional
Sobolev spaces. In this way, our results unify and generalize, to
the best of our knowledge, all previous results in the literature
and explain previous results for indicator functions by their Sobolev
regularity. In particular, arguing by Sobolev regularity instead of Hölder regularity yields approximations for occupation times and local times, if they exist.

The central idea is to expand the $L^{2}(\P)$-norm of $\Gamma_{T}(f)-\widehat{\Gamma}_{T,n}(f)$
in terms of the bivariate distributions $(X_{t},X_{t'})$, $0\leq t,t'\leq T$,
and to derive upper bounds in terms of either their Lebesgue densities
or their characteristic functions. This approach is therefore generic
and not restricted to Markov processes, and covers fractional Brownian
motion and other non-Markovian processes.

For the $L^{2}(\P)$-error, lower bounds can be derived by the conditional
expectation of $\Gamma_{T}(f)$ with respect to the data. For Brownian
motion and functions with fractional Sobolev regularity, this idea
is used to prove that the upper bounds are sharp with respect to $\Delta_{n}$
and $T$. In particular, no other quadrature rule can achieve a faster
rate of convergence than the Riemann estimator uniformly over the
considered function class. Deriving similar upper and lower bounds
for strong $L^{p}(\P)$-approximations and $p>2$ for processes different
from Brownian motion is a challenging problem left for future research. 

Let us shortly review related findings in the literature. Central
limit theorems for $\Gamma_{T}(f)-\widehat{\Gamma}_{T,n}(f)$ were
studied in the case of semimartingales in Chapter 6 of \citep{jacod2011discretization}
for $f\in C^{2}(\R^{d})$ and by \citep{Altmeyer:2019wd} for weakly
differentiable functions. The weak error $\E[\Gamma_{T}(f)-\widehat{\Gamma}_{T,n}(f)]$
was considered by \citep{Gobet2008,Kohatsu-Higa2014} for bounded
$f$, and here higher regularity of $f$ does not improve the result.
Using different techniques, \citep{Ganychenko2015a,Ganychenko2015,Ganychenko2014,Kohatsu-Higa2014}
study the $L^{p}(\P)$-error for Hölder functions and Markov processes
satisfying heat kernel bounds, in particular, scalar diffusions. \citep{Ngo2011,Kohatsu-Higa2014,Jacod1998,ivanovs2020optimal}
approximate occupation times $\Gamma_{T}(\I_{[a,b]})$, $a,b\in\R$
and local times for scalar diffusions, while \citep{Anonymous:2019wc}
estimate derivatives of local times for fractional Brownian motion.
Surprisingly, for occupation times and Brownian motion the rate of
convergence corresponds to the one obtained for occupation time functionals
and $1/2$-Hölder-continuous functions, which cannot be explained
by the specific analysis for indicator functions. For stationary diffusion
processes with infinitesimal generator in divergence form this is
achieved by \citep{Altmeyer:2017dza}, who consider fractional Sobolev
spaces. Since the proof relies heavily on stationarity and semigroup
theory, it is not clear how this can be generalized. We want to emphasize that none of the results above applies to general occupation time functionals and non-Markovian processes. 

This paper is organized as follows. Section \ref{sec:Upper-bounds}
derives general upper bounds for bounded or square integrable functions
$f$. In Section \ref{sec:applications}, several concrete processes $X$ are studied, namely
fractional Brownian motion, Markov processes and processes with independent
increments. This does not cover all possible examples, by far, but
hints at how to derive similar results for other processes. The reader
interested in scalar fractional Brownian motion only may skip to Theorem
\ref{thm:fracBM_d_1} below. The approximation of occupation and local
times is proven only for fractional Brownian motion (Corollaries \ref{cor:occupationTime},
\ref{cor:localTime}), but the proofs are generic and apply to all
other examples. Finally, Section \ref{sec:Lower-bounds} shows that
the upper bounds are sharp for Brownian motion. 

Proofs are deferred to the appendix. Let us introduce some notation.
In the following, $C$ always denotes a positive absolute constant,
which may change from line to line. We write $a\lesssim b$ for $a\leq Cb$
and set $\lfloor t\rfloor_{\Delta_{n}}:=\lfloor t/\Delta_{n}\rfloor\Delta_{n}$
for $t\geq0$, $\norm{\cdot}_{\infty}$ is the usual sup-norm on $\R^{d}$
and $\norm{\cdot}_{L^{2}}$ the $L^{2}$-norm. For $0\leq s\leq1$
denote by $C^{s}(\R^{d})$ and $H^{s}(\R^{d})$ the Hölder and fractional
Sobolev spaces of order $s$, that is, 
\begin{align*}
C^{s}(\R^{d}) & =\left\{ f:\norm f_{\infty}+\norm f_{C^{s}}<\infty\right\} ,\quad\lVert f\rVert_{C^{s}}=\sup_{x\neq y}\frac{|f(x)-f(y)|}{|x-y|^{s}},\\
H^{s}(\R^{d}) & =\left\{ f:\norm f_{L^{2}}+\norm f_{H^{s}}<\infty\right\} ,\quad\norm f_{H^{s}}^{2}=\int_{\R^{d}}\left|\c Ff\left(u\right)\right|^{2}\left|u\right|^{2s}du,
\end{align*}
where $\c Ff$ is the Fourier transform of $f$, which for $f\in L^{1}(\R^{d})\cap L^{2}(\R^{d})$
is given by $\c Ff(u)=\int_{\R^{d}}f(x)e^{i\sc ux}dx$, $u\in\R^{d}$.

\section{General $L^{2}(\protect\P)$-upper bounds\label{sec:Upper-bounds}}

Let $X=(X_{t})_{0\leq t\leq T}$ be a càdlàg process on a filtered
probability space $(\Omega,\c F,(\c F)_{0\leq t\leq T},\P)$ and let
$f:\R^{d}\rightarrow\R$ be measurable. We want to explicitly allow
discontinuous functions $f$. For this we make the standing assumption
that the compositions $f(X_{t})$ are well-defined random variables.
In particular, $f$ bounded or $f\in L^{2}(\R^{d})$ means that $f$
is a function and not only an equivalence class. On the other hand,
if the distribution of $X_{t}$ has a Lebesgue density and $f=\tilde{f}$
almost everywhere, then $f(X_{t})=\tilde{f}(X_{t})$ $\P$-almost
surely.

Let $(u,v)\mapsto\varphi(u,t;v,t')=\E[e^{i\sc u{X_{t}}+i\sc v{X_{t'}}}]$
denote the characteristic function of the bivariate random variable
$(X_{t},X_{t'})\in\R^{2d}$, and $0\leq t,t'\leq T$, for $t\neq t'$
let $(x,y)\mapsto p(x,t;y,t')$ denote its Lebesgue density, if it
exists. We always assume that $(x,t,y,t')\mapsto p(x,t;y,t')$ and
$(u,t,v,t')\mapsto\varphi(u,t;v,t')$ are jointly measurable. In order
to handle singularities for the distribution of $(X_{t},X_{t'})$
at $t=t'$ and near $t=0$ (e.g., when $X_{0}$ is constant), introduce
the set
\begin{align}
\mathbb{T}_{T,n} & =\{(t,t')\in[0,T]^{2}:t\geq\Delta_{n},t'>t+2\Delta_{n}\}.\label{eq:T_n}
\end{align}
After these preliminaries we obtain the following general upper bounds. 
\begin{prop}
\label{prop:upperBound_density}Assume that the distributions of the
bivariate random variables $(X_{t},X_{t'})$ have Lebesgue densities
$p(\cdot,t;\cdot,t')$ for all $t\neq t'$, $t,t'>0$. Then the following
holds for bounded $f$:
\begin{enumerate}
\item If $t'\mapsto p(x,t;y,t')$ is differentiable for all $x,y\in\R^{d}$,
$0<t<t'<T$ with $\partial_{t'}p(x,\cdot;y,\cdot)\in L^{1}(\mathbb{T}_{T,n})$,
then
\begin{align*}
\norm{\Gamma_{T}(f)-\widehat{\Gamma}_{T,n}(f)}_{L^{2}(\P)}^{2} & \leq C\int_{\R^{2d}}(f(x)-f(y))^{2}(w_{T,n}^{(x,y)}+v_{T,n}^{(x,y)})d(x,y),\\
\text{with}\,\,\,w_{T,n}^{(x,y)}=\Delta_{n} & \int_{\Delta_{n}}^{T}p(x,\lfloor t\rfloor_{\Delta_{n}};y,t)dt+T\int_{0}^{\Delta_{n}}p(x,\lfloor t\rfloor_{\Delta_{n}};y,t)dt,\\
v_{T,n}^{(x,y)}=\Delta_{n} & \int_{\mathbb{T}_{T,n}}|\partial_{t'}p(x,t;y,t')-\partial_{t'}p(x,\lfloor t\rfloor_{\Delta_{n}};y,t')|d(t,t').
\end{align*}
\item If also $t\mapsto\partial_{t'}p(x,t;y,t')$ is differentiable for
$t<t'$ with $\partial_{tt'}^{2}p(x,\cdot;y,\cdot)\in L^{1}(\mathbb{T}_{T,n})$,
then the upper bound in (i) holds with $v_{T,n}^{(x,y)}$ replaced
by
\begin{align*}
\tilde{v}_{T,n}^{(x,y)}= & \Delta_{n}^{2}\int_{\mathbb{T}_{T,n}}|\partial_{rr'}^{2}p(x,t;y,t')|d(t,t').
\end{align*}
\end{enumerate}
\end{prop}

The proof of this proposition is inspired by Theorem 1 of \citep{Ganychenko2015a},
but applies also to non-Markovian processes. In order to formulate
a similar result with respect to the characteristic functions, we
include an additional smoothing by an independent random variable
$\xi$.
\begin{prop}
\label{prop:upperBound_FT} Let $f\in L^{2}(\R^{d})$ and let $\xi$
be an $\R^{d}$-valued random variable, independent of $X$ with bounded
Lebesgue density $\mu$.
\begin{enumerate}
\item If $t'\mapsto\varphi(u,t;v,t')$ is differentiable for all $u,v\in\R^{d}$,
$0<t<t'<T$ with $t'\mapsto\partial_{t'}\varphi(u,t;v,t')\in L^{1}(\mathbb{T}_{T,n})$,
then
\begin{align*}
 & \norm{\Gamma_{T}(f(\cdot+\xi))-\widehat{\Gamma}_{T,n}(f(\cdot+\xi))}_{L^{2}(\P)}^{2}\leq C\norm{\mu}_{\infty}\int_{\R^{d}}|\c Ff(u)|^{2}(w_{T,n}^{(u)}+v_{T,n}^{(u)})du,\\
 & \qquad\text{with}\,\,\,w_{T,n}^{(u)}=\Delta_{n}\int_{\Delta_{n}}^{T}|g(u,t)|dt+T\int_{0}^{\Delta_{n}}|g(u,t)|dt,\\
 & \qquad\qquad v_{T,n}^{(u)}=\Delta_{n}\int_{\mathbb{T}_{T,n}}|\partial_{t'}\varphi(u,t;-u,t')-\partial_{t'}\varphi(u,\lfloor t\rfloor_{\Delta_{n}};-u,t')|d(t,t'),
\end{align*}
where $g(u,t)=2-2\text{Re}\varphi(u,\lfloor t\rfloor_{\Delta_{n}};-u,t)$.
\item If also $t\mapsto\partial_{t'}\varphi(u,t;v,t')$ is differentiable
for $t<t'$ with $t\mapsto\partial_{tt'}^{2}\varphi(u,t;v,t')\in L^{1}(\mathbb{T}_{T,n})$,
then the upper bound in (i) holds with $v_{T,n}^{(u)}$ replaced by
\begin{align*}
\tilde{v}_{T,n}^{(u)} & =\Delta_{n}^{2}\int_{\mathbb{T}_{T,n}}|\partial_{tt'}^{2}\varphi(u,t;-u,t')|d(t,t').
\end{align*}
\end{enumerate}
\end{prop}

Due to these two results the $L^{2}(\P)$-error is controlled by upper
bounding
\begin{align*}
\int_{\R^{2d}}(f(x)-f(y))^{2}\rho(x,y)d(x,y), & \int_{\R^{d}}|\c Ff(u)|^{2}\bar{\rho}(u)du,
\end{align*}
with respect to certain weight functions $\rho$, $\bar{\rho}$. This
naturally suggests to consider $f\in C^{s}(\R^{d}),H^{s}(\R^{d})$
and to prove that the weights decay sufficiently fast. For example, if $f\in H^{s}(\R^{d})$,
then 
\begin{equation}
\int_{\R^{d}}|\c Ff(u)|^{2}\bar{\rho}(u)du\leq\norm f_{H^{s}}^{2}\sup_{u}|u|^{-2s}\bar{\rho}(u).\label{eq:fourierUpper}
\end{equation}
If also $\rho(x,y)=\rho(y-x)$, then the Plancherel theorem shows
\begin{align}
 & \int_{\R^{2d}}(f(x)-f(y))^{2}\rho(y-x)d(x,y)=(2\pi)^{-2d}\int_{\R^{2d}}|\c Ff(u)(1-e^{-i\sc uy})|^{2}\rho(y)d(u,y)\nonumber \\
 & \quad\leq(2\pi)^{-2d}\norm f_{H^{s}}^{2}\int_{\R^{d}}|y|^{2s}\rho(y)dy,\label{eq:densityFourierUpper}
\end{align}
and so we are left with with studying the weights $\rho$, $\bar{\rho}$.
\begin{rem}
\label{rem:xi}%
\mbox{%
}
\begin{enumerate}
\item The regularization with $\xi$ in Proposition \ref{prop:upperBound_FT}
has also been used in Theorem 2 of \citep{Altmeyer:2019wd}. The independence
of $\xi$ allows for using $L^{2}$-arguments in the proofs, cf. inequality
(\ref{eq:boundedDensity}) below. We will use in Theorem \ref{thm:fracBM_d_1}
below a different argument involving Fourier transforms, that demonstrates
how we can argue, in principle, without this additional smoothing.
Moreover, no smoothing is necessary for deriving the lower bounds
in Section \ref{sec:Lower-bounds}.
\item The upper bounds in Proposition \ref{prop:upperBound_FT} depend only
on $\varphi(u,t;-u,t')=\E[e^{i\sc u{X_{t}-X_{t'}}}]$, that is, on
increments $X_{t'}-X_{t}$.
\end{enumerate}
\end{rem}

\section{Application to examples}\label{sec:applications}

\subsection{\label{subsec:Fractional-Brownian-motion}Fractional Brownian motion}

Let $X$ be a fractional Brownian motion in $\R^{d}$ with Hurst index
$0<H<1$. The $d$ component processes $(X_{t}^{(m)})_{0\leq t\leq T}$
for $m=1,\dots,d$ are independent centered Gaussian processes with
covariance function 
\begin{equation}
c(t,t')=\E[X_{t}^{(m)}X_{t'}^{(m)}]=\frac{1}{2}((t')^{2H}+t^{2H}-|t'-t|^{2H}),\,\,\,0\leq t,t'\leq T.\label{eq:fracBMcovFun}
\end{equation}
For $H=1/2$, $X$ is a Brownian motion. For $H\neq1/2$, fractional
Brownian motion is an important example of a non-Markovian process,
which is also not a semimartingale. 

Both the densities and the characteristic functions of $(X_{t},X_{t'})$,
$0<t<t'\leq T$, are explicit by Gaussianity, but it is much easier
to upper bound the time derivatives of the latter one. In the setting
of Proposition \ref{prop:upperBound_FT} we have:
\begin{thm}
\label{thm:fractionalBM} Let $X$ be a $d$-dimensional fractional
Brownian motion with $0<H<1$, and let $X_{0}=\xi$ be as in Proposition
\ref{prop:upperBound_FT}. If $f\in H^{s}(\R^{d})$, $0\leq s\leq\min(1,1/(2H))$,
then
\begin{align*}
\norm{\Gamma_{T}(f)-\widehat{\Gamma}_{T,n}(f)}_{L^{2}(\P)} & \leq C\norm{\mu}_{\infty}^{1/2}\norm f_{H^{s}}T^{1/2}\Delta_{n}^{1/2+sH}.
\end{align*}
\end{thm}

For $H=1/2$, this extends Theorem 3.9 of \citep{Altmeyer:2017dza}
to general dimension $d$. We see that the rate $\Delta_{n}^{1/2+sH}$
improves for more regular $f$, but never exceeds $\Delta_{n}$. In
Theorem \ref{thm:lowerBound_Hs} we show for $H=1/2$ that it is optimal.

In the scalar case, we demonstrate now how the random variable $\xi$
in Proposition \ref{prop:upperBound_FT} can be avoided. The proof
is not based on Proposition \ref{prop:upperBound_FT} and replaces
the role of the density $\mu$ by the characteristic function $\varphi(u,t;v,t')$
to achieve sufficient integrability, cf. Remark \ref{rem:xi}.
\begin{thm}
\label{thm:fracBM_d_1} Let $X$ be a scalar fractional Brownian motion
with $0<H\leq1/2$. If $f\in H^{s}(\R)$, $0\leq s\leq1$, is bounded,
then
\begin{align*}
 & \norm{\Gamma_{T}(f)-\widehat{\Gamma}_{T,n}(f)}_{L^{2}(\P)}\leq C(\norm f_{\infty}\Delta_{n}+\norm f_{H^{s}}^{1/2}\max(\norm f_{H^{s}}^{1/2},\norm f_{L^{2}}^{1/2})\\
 & \qquad\max(1,T^{1/4})T^{1/4-H/2}\Delta_{n}^{1/2+sH}).
\end{align*}
\end{thm}

While the rate $\Delta_{n}^{1/2+sH}$ remains unchanged, the factor
$T^{1/2}$ from Theorem \ref{thm:fractionalBM} changes to $T^{1/2-H/2}$
for $T\geq1$. The proof of the lower bound in Theorem \ref{thm:lowerBound_Hs}
reveals that this is necessary for $H=1/2$ due to the singularity
of the distribution of $X_{t}$ near $t=0$. Moreover, for $H>1/2$,
the proof of Theorem \ref{thm:fracBM_d_1} yields a slower rate than
$\Delta_{n}^{1/2+sH}$.

Our interest in deriving convergence rates for $f$ with fractional
Sobolev regularity is the possibility to approximate occupation times
and local times of a scalar process $X$. Such rates will now be derived
as corollaries from Theorem \ref{thm:fracBM_d_1}. With respect to
the occupation time of a set $[a,b]$, $-\infty<a<b<\infty$, note
that $f=\I_{[a,b]}\in H^{1/2-}(\R)$. Theorem \ref{thm:fracBM_d_1}
therefore suggests the rate $\Delta_{n}^{1/2+H/2-}$, but by interpolation
we obtain even an exact result.
\begin{cor}
\label{cor:occupationTime}Let $X$ be a scalar fractional Brownian
motion with $0<H<1$ and let $-\infty\leq a<b\leq\infty$. If $X_{0}=\xi$
is as in Proposition \ref{prop:upperBound_FT}, then the Riemann estimator
for the occupation time of the set $[a,b]$ satisfies
\begin{align*}
\norm{\Gamma_{T}(\I_{[a,b]})-\widehat{\Gamma}_{T,n}(\I_{[a,b]})}_{L^{2}(\P)} & \leq C\norm{\mu}_{\infty}^{1/2}T^{1/2}\Delta_{n}^{1/2+H/2},
\end{align*}
with a constant $C$ independent of $a,b$. If $0<H\leq1/2$ and $-\infty<a<b<\infty$,
then also with $X_{0}=0$
\begin{align*}
\norm{\Gamma_{T}(\I_{[a,b]})-\widehat{\Gamma}_{T,n}(\I_{[a,b]})}_{L^{2}(\P)} & \leq C_{a,b}\max(1,T^{1/4})T^{1/4-H/2}\Delta_{n}^{1/2+H/2},
\end{align*}
but this time the constant $C_{a,b}$ depends on $a,b$.
\end{cor}

On the other hand, let $(L_{T}(a))_{a\in\R}$ be the family of \emph{local
times }of $X$ until $T$, cf. Chapter 5 of \citep{nualart1995malliavin}.
Formally, $L_{T}(a)=\Gamma_{T}(\delta_{a})$, where $\delta_{a}$
is the Dirac delta function. If we use for $H^{s}(\R)$ the norm $\vertiii{f}{}_{H^{s}}=(\int_{\R}|\c Ff(u)|^{2}(1+|u|^{2})^{s}du)^{1/2}$,
then this allows for defining Sobolev spaces of negative regularity
$s<0$, which is a space of distributions. In particular, the Fourier
transform (in the sense of distributions) of $\delta_{a}$ is $\c F\delta_{a}(u)=e^{-iua}$,
implying $\delta_{a}\in H^{s}(\R)$ for $s<-1/2$, and so Theorem
\ref{thm:fracBM_d_1} formally suggests the rate $\D^{1/2-H/2-}$. For simplicity,
we only consider $0<H\leq1/2$, the case $1/2<H<1$ follows as
in Corollary \ref{cor:occupationTime} assuming $X_{0}=\xi$.
\begin{cor}
\label{cor:localTime} Let $X$ be a scalar fractional Brownian motion
with $0<H\leq1/2$ and set $f_{a,n}(x)=(2\Delta_{n}^{H})^{-1}\I_{[a-\Delta_{n}^{H},a+\Delta_{n}^{H}]}(x)$
for $x,a\in\R$. Then for any $\epsilon>0$
\begin{align*}
\norm{L_{T}(a)-\widehat{\Gamma}_{T,n}(f_{a,n})}_{L^{2}(\P)} & \leq C\max(1,T^{1/2})T^{1/4-H/2}\Delta_{n}^{1/2-H/2-\epsilon}.
\end{align*}
\end{cor}

For $H=1/2$, the results of Corollaries \ref{cor:occupationTime},
\ref{cor:localTime} have been obtained before (\citep{Ngo2011},
\citep{Kohatsu-Higa2014}), but are new for general $H$. The loss
in the rate $\Delta_{n}^{1/2-H/2-\epsilon}$ in Corollary \ref{cor:localTime}
for arbitrarily small $\epsilon>0$ is due to a suboptimal moment
bound for local times. In case of a Brownian motion, we can instead
use \citep{marcus2008lp} to remove $\epsilon$, and recover the rate
$\D^{1/4}$ from \citep{Jacod1998}, improving on Theorem 2.6 of \citep{Kohatsu-Higa2014}.

\subsection{\label{subsec:Markov-processes}Markov processes with heat kernel
bounds}

Let $X$ be a continuous-time Markov process on $\R^{d}$ with transition
densities $p_{t,t'}$, $0\leq t<t'\leq T$, such that 
\[
\E[f(X_{t'})|X_{t}=x]=\int_{\R^{d}}f(y)p_{t,t'}(x,y)dy,\,\,\,x\in\R^{d},
\]
for continuous and bounded $f$. The distribution of $(X_{t},X_{t'})$,
$t<t'$, conditional on $X_{0}=x_{0}\in\R^{d}$, has the density $p(x,t;y,t';x_{0})=p_{0,t}(x_{0},x)p_{t,t'}(x,y)$.
Suppose the following:

\begin{assumption}\label{ass:A} $X$ is a $d$-dimensional Markov
process such that $(x,t,y,t')\mapsto p_{t,t'}(x,y)$ is jointly measurable
and $t'\mapsto p_{t,t'}(x,y)$ is continuously differentiable for
all $x,y\in\R^{d}$, $0<t<t'<T$, and such that the heat kernel bounds
\begin{align*}
p_{t,t'}(x,y)\leq Cq_{t'-t}(y-x),\,\, & |\partial_{t'}p_{t,t'}(x,y)|\leq C\frac{q_{t'-t}(y-x)}{t'-t},
\end{align*}
are satisfied, where the $(q_{t})_{0<t\leq T}$ are probability densities
and $(t,x)\mapsto q_{t}(x)$ is jointly measurable. \end{assumption}

\begin{assumption}\label{ass:B} In addition to Assumption \ref{ass:A},
$t\mapsto\partial_{t'}p_{t,t'}(x,y)$ is continuously differentiable
for all $x,y\in\R^{d}$, $0<t<t'<T$, $\int_{\R^{d}}|x|^{\alpha'}q_{t}(x)dx\leq Ct^{\alpha'/\alpha}$
for some $0<\alpha\leq2$ and all $0<\alpha'\leq\alpha$, and
\[
|\partial_{tt'}^{2}p_{t,t'}(x,y)|\leq C\frac{q_{t'-t}(y-x)}{(t'-t)^{2}}.
\]
\end{assumption}

Such heat kernel bounds are satisfied for elliptic diffusion processes
with sufficiently regular coefficients. In this case the transition
densities satisfy the Kolmogorov forward equation
\[
\partial_{t'}p_{t,t'}(x,\cdot)=L^{*}p_{t,t'}(x,\cdot),\,\,t<t',\,x\in\R^{d},
\]
where $L^{*}$ is the adjoint of the infinitesimal generator of $X$,
cf. Chapter 5.7 of \citep{Karatzas1991}. Upper bounds on $\partial_{t'}p_{t,t'}$,
$\partial_{tt'}^{2}p_{t,t'}$ follow therefore from bounds on the
partial derivatives of $(x,y)\mapsto p_{t,t'}(x,y)$ with $q_{t}(x)=Ct^{-d/2}e^{-C|xt^{-1/2}|^{2}}$,
$\alpha=2$, cf. Theorem 9.4.2 of \citep{Friedman:2013uk}. Important
examples for Markov processes satisfying Assumptions \ref{ass:A}
and \ref{ass:B} are Lévy driven SDEs ( \citep{Knopova:2018ex,Kuhn:2019cw}),
in particular, $\alpha$-stable processes. For slightly more general
heat kernel bounds see \citep{Ganychenko2015}.

Plugging the heat kernel bounds into the abstract bounds of Proposition
\ref{prop:upperBound_density} yields the following. We write $\P_{x_{0}}$
to indicate the initial value $X_{0}=x_{0}$.
\begin{thm}
\label{thm:density_L2_Hs} Let $X_{0}$ have a bounded Lebesgue density
$\mu$. Under Assumption \ref{ass:A} we have for all $f\in L^{2}(\R^{d})$
\[
\norm{\Gamma_{T}(f)-\widehat{\Gamma}_{T,n}(f)}_{L^{2}\left(\P\right)}\leq C\norm{\mu}_{\infty}^{1/2}\norm f_{L^{2}}T^{1/2}\Delta_{n}^{1/2}(\log n)^{1/2},
\]
while under Assumption \ref{ass:B} for $f\in H^{s}(\R^{d})$, $0\leq s\leq\alpha/2$,
\[
\norm{\Gamma_{T}(f)-\widehat{\Gamma}_{T,n}(f)}_{L^{2}(\P)}\leq C\norm{\mu}_{\infty}^{1/2}\norm f_{H^{s}}\begin{cases}
T^{1/2}\Delta_{n}^{1/2+s/\alpha}, & s<\alpha/2,\\
T^{1/2}\Delta_{n}(\log n)^{1/2}, & s=\alpha/2.
\end{cases}
\]
\end{thm}

\begin{thm}
\label{thm:density_bdd_Holder}Under Assumption \ref{ass:A} we have
for bounded $f$
\[
\norm{\Gamma_{T}(f)-\widehat{\Gamma}_{T,n}(f)}_{L^{2}(\P_{x_{0}})}\leq C\norm f_{\infty}T^{1/2}\Delta_{n}^{1/2}(\log n)^{1/2},
\]
while under Assumption \ref{ass:B} for $f\in C^{s}(\R^{d}),0\leq s\leq\alpha/2$,
\[
\norm{\Gamma_{T}(f)-\widehat{\Gamma}_{T,n}(f)}_{L^{2}(\P_{x_{0}})}\leq C\norm f_{C^{s}}\begin{cases}
T^{1/2}\Delta_{n}^{1/2+s/\alpha}, & s<\alpha/2,\\
T^{1/2}\Delta_{n}(\log n)^{1/2}, & s=\alpha/2.
\end{cases}
\]
\end{thm}

These two results unify and generalize previously obtained results
for the strong $L^{2}(\P)$-error for Markov diffusion processes in
arbitrary dimension $d$, up to log-terms. These are Theorem 2.14
of \citep{Ganychenko2014}, Theorem 2.1 \citep{Ganychenko2015} for
bounded $f$, \citep{Ganychenko2015a}, Theorem 2.3 of \citep{Kohatsu-Higa2014}
for $f\in C^{s}(\R^{d})$ and Theorem 3.7 of \citep{Altmeyer:2017dza}
for $f\in H^{s}(\R^{d})$. The additional $\log$-terms are not optimal,
in general. They are not present for $f\in H^{s}(\R^{d})$, for example,
when $X$ is a Brownian motion (Theorem \ref{thm:fractionalBM}) or
a stationary reversible diffusion process (Theorem 3.7 of \citep{Altmeyer:2017dza}).

Results for occupation or local times can be obtained exactly as in
Corollaries \ref{cor:occupationTime}, \ref{cor:localTime} by interpolation
in Theorem \ref{thm:density_L2_Hs}, and as long as moment bounds
for the local times exist. We want to emphasize that this interpolation does not work using  Theorem \ref{thm:density_bdd_Holder}, since indicator functions are not well approximated in Hölder spaces.

\begin{rem}
There are different ways to relax the assumption on $X_{0}$ in Theorem
\ref{thm:density_L2_Hs}.
\begin{enumerate}
\item If we are estimating $\int_{T_{0}}^{T}f(X_{t})dt$ for $T_{0}>0$
using the corresponding Riemann estimator, then Theorem \ref{thm:density_L2_Hs}
remains valid by the Markov property, if $X_{0}$ is replaced by $X_{T_{0}}$
and $q$ from Assumption \ref{ass:A} is bounded.
\item It is enough to upper bound $\sup_{x\in\R^{d}}q_{t}(x)$ in the proof
(cf. Equations (\ref{eq:heat_bound_1}) to (\ref{eq:heat_bound_3})
and (\ref{eq:heat_bound_1-1}) to (\ref{eq:heat_bound_3-1})). Since
$q_{t}$ typically has a singularity near $t=0$, this will yield
a slower rate.
\item If $X_{0}\in L^{2}(\R^{d})$ is not bounded, then we obtain again
slower rates, cf. \citep{Neuenkirch:2019uo}.
\end{enumerate}
\end{rem}

\begin{rem}
The assumption on $f$ being bounded in Theorem \ref{thm:density_bdd_Holder}
can be relaxed by considering weighted Hölder norms, cf. \citep{Ganychenko2015,Ganychenko2014}.
\end{rem}

\subsection{\label{subsec:Additive-processes}Processes with independent increments}

Let $X$ be an additive process on $\R^{d}$ with local characteristics
$(\sigma_{t},F_{t},b_{t})_{t\geq0}$, where $t\mapsto\sigma_{t}$
is a continuous $\R^{d\times d}$-valued function, $t\mapsto b_{t}$
is a locally integrable $\R^{d}$-valued function and $(F_{t})_{t\geq0}$
is a family of positive measures on $\R^{d}$ with $F_{t}(\{0\})=0$
and $\sup_{0\leq t\leq T}\{\int(|x|^{2}\wedge1)dF_{t}(x)\}<\infty$,
cf. Chapter 14 of \citep{tankov2003financial}.

$X$ is an inhomogeneous Markov process with independent increments,
in particular every Lévy process is an additive process. We can therefore
apply the results from Section \ref{subsec:Markov-processes} as soon
as heat kernel bounds are available. In general, however, it is rather
difficult to compute or even upper bound the marginal densities. When
$\sigma_{t}\sigma_{t}^{\top}$ is not invertible at some $t$, the
densities might not even exist. On the other hand, the characteristic
functions are known explicitly by the Lévy-Khintchine formula, cf.
Theorem 14.1 of \citep{tankov2003financial}. The characteristic function
of $(X_{t},X_{t'})$, $0<t<t'$, is $\varphi(u,t;v,t')=e^{\Psi_{t,t'}(v)+\Psi_{0,t}(u+v)}$,
$u,v\in\R^{d}$, with characteristic exponents $\Psi_{t,t'}(u)$ equal
to
\[
i\int_{t}^{t'}\sc u{b_{s}}ds-\frac{1}{2}\int_{t}^{t'}|\sigma_{s}^{\top}u|^{2}ds+\int_{t}^{t'}\int_{\R^{d}}(e^{i\sc ux}-1-i\sc ux\I_{\left\{ |x|\leq1\right\} })dF_{s}(x)ds.
\]
For concrete bounds suppose the following:

\begin{assumption}\label{assu:II} $X$ is a $d$-dimensional additive
process with characteristic exponents satisfying for $0\leq t<t'\leq T$,
$u\in\R^{d}$
\begin{align*}
|e^{\Psi_{t,t'}(u)}| & \leq Ce^{-C|u|^{\alpha}(t'-t)},\,\,|\Psi_{t,t'}(u)|\leq C\max(1,|u|^{\alpha^{*}})|t'-t|,
\end{align*}
where $0\leq\alpha\leq2$, $\beta\geq0$ are such that $0\leq\alpha(1+\beta)\leq2$,
$\alpha^{*}=\max(1,\alpha(1+\beta))$.\end{assumption}

This assumption holds, for example, if $X$ is a generalized $\alpha$-stable
process with $\alpha^{*}=\alpha$ or with time varying stability index
$t\mapsto\alpha(1+\beta_{t})$, $0\leq\beta_{t}\leq\beta$. On the
other hand, if $\sigma_{t}\sigma_{t}^{\top}$ is non-degenerate for
all $t$, then $\alpha^{*}=\alpha=2$ (use Equation 8.9 of \citet{sato1999levy}).
\begin{thm}
\label{thm:additiveProcess} Grant Assumption \ref{assu:II} and let
$X_{0}$ have a bounded Lebesgue density $\mu$. Then:
\begin{enumerate}
\item If $f\in H^{s}(\R^{d})$, $0\leq s\leq\alpha^{*}/2$, then 
\[
\norm{\Gamma_{T}(f)-\widehat{\Gamma}_{T,n}(f)}_{L^{2}(\P)}\leq C\norm{\mu}_{\infty}^{1/2}\max(\norm f_{L^{2}},\norm f_{H^{s}})T^{1/2}\Delta_{n}^{1/2+s/\alpha^{*}}.
\]
\item If $\alpha^{*}=0$, then for $f\in L^{2}(\R^{d})$ the upper bound
is $C\norm{\mu}_{\infty}^{1/2}\norm f_{L^{2}}(T^{1/2}+T)\Delta_{n}$.
If $-C(t'-t)\leq\Psi_{t,t'}(u)\leq0$ for $t'>t$, then $T^{1/2}+T$
can be replaced by $T^{1/2}$.
\end{enumerate}
\end{thm}

If $\alpha^{*}=\alpha$, then the rate in (i) is $\Delta_{n}^{1/2+s/\alpha}$
as in Theorem \ref{thm:density_L2_Hs}, but without the $(\log n)^{1/2}$-term.
$\alpha^{*}=0$ holds for a compound Poisson process. In this case,
$\Psi_{t,t'}(u)=|t'-t|\int_{\R^{d}}(e^{i\sc ux}-1)dF(x)$ for a finite
measure $F$, and so $\Psi_{t,t'}(u)$ is bounded. The improved bound
in (ii) applies, if $F$ is symmetric. For stationary $X$ this has
been shown also in Section 3.1 of \citep{Altmeyer:2017dza}. As above,
we can argue as in Corollaries \ref{cor:occupationTime} and \ref{cor:localTime}
to obtain estimates for occupation times and local times of $X$,
as long as the latter exist.

\section{\label{sec:Lower-bounds} Sharpness of upper bounds for Brownian
motion}

Recall from Theorems \ref{thm:fractionalBM} and \ref{thm:fracBM_d_1}
that the upper bound for the $L^{2}(\P)$-error with respect to the
Riemann estimator $\widehat{\Gamma}_{T,n}(f)$ is of order $\Delta_{n}^{(1+s)/2}$
as $n\rightarrow\infty$ for $f\in H^{s}(\R^{d})$ and when $X$ is
a Brownian motion. We prove now that this rate is sharp up to a log-factor
uniformly in $H^{s}(\R^{d})$.

The key idea of the following proof, which extends Theorem 5 of \citep{Altmeyer:2019wd},
is that the minimal $L^{2}(\P)$-error is achieved by the conditional
expectation, 
\[
\norm{\Gamma_{T}(f)-\widehat{\Gamma}}_{L^{2}(\P)}\geq\norm{\Gamma_{T}(f)-\E[\Gamma_{T}(f)|\mathcal{G}_{n}]}_{L^{2}(\P)},
\]
where $\widehat{\Gamma}$ is any square integrable estimator of $\Gamma_{T}(f)$,
being measurable with respect to the sigma field $\c G_{n}=\sigma(X_{t_{k}}:k\in\{0,\dots,n\})$.
When $X$ is a Brownian motion, the conditional expectation can be
computed explicitly. The only other explicit lower bound in the literature
is Proposition 2.3 of \citep{Ngo2011} for $f=\I_{[0,\infty)}$ and
a scalar Brownian motion.
\begin{thm}
\label{thm:lowerBound_Hs} Let $X$ be a $d$-dimensional Brownian
motion and define for $0\leq s\leq1$
\[
f^{*}=\c F^{-1}g,\quad\text{with }g(u)=\frac{1}{(1+|u|)^{s+d/2}\log(1+|u|)},u\in\R^{d}.
\]
Then $f^{*}\in H^{s}(\R^{d})$ and we have the asymptotic lower bound
\[
\liminf_{n\rightarrow\infty}\inf_{\hat{\Gamma}}\left(\log(\Delta_{n}^{-1/2})\Delta_{n}^{-(1+s)/2}\norm{\Gamma_{T}(f^{*})-\hat{\Gamma}}_{L^{2}(\P)}\right)\ge CT^{1/4}>0,
\]
where $\hat{\Gamma}$ is any square integrable estimator for $\Gamma_{T}(f^{*})$
based on $(X_{t_{k}})_{k\in\{0,\dots,n\}}$. In particular, the rate
$\Delta_{n}^{(1+s)/2}$, achieved by the Riemann estimator $\hat{\Gamma}_{T,n}(f)$,
is sharp up to a log-factor uniformly for $f\in H^{s}(\R^{d})$.
\end{thm}

If $X_0$ is not constant, then the time dependency in the error necessarily increases. 
\begin{cor}
\label{cor:lowerBound}Consider the setting of Theorem \ref{thm:lowerBound_Hs}.
If $X_{0}$ is independent of $X-X_0$ and has a Lebesgue density $\mu$ with nonnegative and Fourier transform $\c F\mu\in L^{1}(\R^{d})$, then
\[
\liminf_{n\rightarrow\infty}\inf_{\hat{\Gamma}}\left(\log(\Delta_{n}^{-1/2})\Delta_{n}^{-(1+s)/2}\norm{\Gamma_{T}(f^{*})-\hat{\Gamma}}_{L^{2}(\P)}\right)\ge CT^{1/2}>0.
\]
\end{cor}

\titleformat{\section}{\Large\bfseries}{\appendixname~\thesection :}{0.5em}{}

\appendix

\section{Proofs}

\subsection{Proof of Proposition \ref{prop:upperBound_density}}
\begin{proof}
(i). Recall the definition of $\mathbb{T}_{T,n}$ from (\ref{eq:T_n})
and set 
\begin{equation}
\tilde{\mathbb{T}}_{n}=\{(t,t')\in[0,T]^{2}:t\geq\Delta_{n},t'>t+4\Delta_{n}\}\label{eq:T_tilde_T_n}
\end{equation}
For $0\leq t,t'\leq T$ let 
\begin{align}
E_{t,t'} & =\E[(f(X_{t})-f(X_{\lfloor t\rfloor_{\Delta_{n}}}))(f(X_{t'})-f(X_{\lfloor t'\rfloor_{\Delta_{n}}}))].\label{eq:E_f}
\end{align}
Using symmetry decompose 
\begin{align}
\norm{\Gamma_{T}(f)-\widehat{\Gamma}_{T,n}(f)}{}_{L^{2}(\P)}^{2} & =\int_{[0,T]^{2}}E_{t,t'}d(t,t')=A_{1}+2A_{2}+2A_{3},\label{eq:A_1_2_3}\\
\text{with } & A_{1}=\int_{[\Delta_{n},T]^{2}}\I_{\{|t-t'|\leq3\Delta_{n}\}}E_{t,t'}d(t,t'),\nonumber \\
 & A_{2}=\int_{\Delta_{n}}^{T}\int_{0}^{\Delta_{n}}E_{t,t'}dt'dt,\,\,A_{3}=\int_{\mathbb{\tilde{T}}_{n}}E_{t,t'}d(t,t').\nonumber 
\end{align}
For the result it is enough to show with
\begin{align*}
w_{T,n}^{(x,y)} & =\Delta_{n}\int_{\Delta_{n}}^{T}p(x,\lfloor t\rfloor_{\Delta_{n}};y,t)dt+T\int_{0}^{\Delta_{n}}p(x,\lfloor t\rfloor_{\Delta_{n}};y,t)dt,\\
v_{T,n}^{(x,y)} & =\Delta_{n}\int_{\mathbb{T}_{T,n}}|\partial_{t'}p(x,t;y,t')-\partial_{t'}p(x,\lfloor t\rfloor_{\Delta_{n}};y,t')|d(t,t'),
\end{align*}
that
\begin{align}
|A_{1}|+|A_{2}| & \lesssim\int_{\R^{2d}}(f(x)-f(y))^{2}w_{T,n}^{(x,y)}d(x,y),\label{eq:A_1_2}\\
|A_{3}| & \lesssim\int_{\R^{2d}}(f(x)-f(y))^{2}v_{T,n}^{(x,y)}d(x,y).\label{eq:A_3}
\end{align}
For the first part a rough argument suffices. Observe that $|E_{t,t'}|\leq\frac{1}{2}E_{t,t}+\frac{1}{2}E_{t',t'}$.
The claim in (\ref{eq:A_1_2}) follows therefore from $|A_{1}|\lesssim\Delta_{n}\int_{\Delta_{n}}^{T}E_{t,t}dt$,
\begin{align}
|A_{2}| & \lesssim\int_{\Delta_{n}}^{T}\int_{0}^{\Delta_{n}}(E_{t,t}+E_{t',t'})dt'dt\leq\Delta_{n}\int_{0}^{T}E_{t,t}dt+T\int_{0}^{\Delta_{n}}E_{t,t}dt,\label{eq:A_2}
\end{align}
and $E_{t,t}=\int_{\R^{2d}}(f(x)-f(y))^{2}p(x,\lfloor t\rfloor_{\Delta_{n}};y,t)d(x,y)$.
With respect to (\ref{eq:A_3}), the regularity assumptions on the
joint densities are crucial. Consider $(t,t')\in\tilde{\mathbb{T}}_{n}$.
Clearly, 
\begin{align*}
 & E_{t,t'}=\int_{\R^{2d}}f(x)f(y)\{p(x,t;y,t')-p(x,t;y,\lfloor t'\rfloor_{\Delta_{n}})\\
 & \qquad-p(x,\lfloor t\rfloor_{\Delta_{n}};y,t')+p(x,\lfloor t\rfloor_{\Delta_{n}};y,\lfloor t'\rfloor_{\Delta_{n}})\}d(x,y).
\end{align*}
Here comes the main insight: If $f(x)$ is replaced in this equality
by $f(y)$, then the $d(x,y)$-integral vanishes. The same holds with
$f(y)$ replaced by $f(x)$. This allows two modifications in the
last display. First, replace $f(x)f(y)$ by $-1/2(f(x)-f(y))^{2}$,
and second, use differentiability of the joint density. Then the last
display reduces to 
\begin{align}
 & -\frac{1}{2}\int_{\R^{2d}}(f(x)-f(y))^{2}\int_{\lfloor t'\rfloor_{\Delta_{n}}}^{t'}\{\partial_{r}p(x,t;y,r)-\partial_{r}p(x,\lfloor t\rfloor_{\Delta_{n}};y,r)\}dr\,d(x,y).\label{eq:ganychenko_trick}
\end{align}
If $(t,t')\in\mathbb{\tilde{T}}_{n}$ and $\lfloor t'\rfloor_{\Delta_{n}}\leq r<t'$,
then $(t,r)\in\mathbb{T}_{T,n}$ and $|t'-r|\leq\Delta_{n}$. Integrating
in the last display over $t,t'\in\tilde{\mathbb{T}}_{n}$, it can
therefore be upper bounded by a double integral over $(t,r)\in\mathbb{T}_{T,n}$,
yielding an additional $\Delta_{n}$. This implies (\ref{eq:A_3}).

(ii). It is enough to prove (\ref{eq:A_3}) with $\tilde{v}_{T,n}^{(x,y)}=\Delta_{n}^{2}\int_{\mathbb{T}_{T,n}}|\partial_{rr'}^{2}p(x,t;y,t')|d(t,t')$
instead of $v_{T,n}^{(x,y)}$. As in (i), $(t,t')\in\mathbb{\tilde{T}}_{n}$
and $\lfloor t\rfloor_{\Delta_{n}}\leq r<t$, $\lfloor t'\rfloor_{\Delta_{n}}\leq r'<t'$
imply $(r,r')\in\mathbb{T}_{T,n}$ and $|t-r|,|t'-r'|\leq\Delta_{n}$.
Since $t\mapsto\partial_{r}p(x,t;y,r)$ is differentiable, the $dr$-integral
in (\ref{eq:ganychenko_trick}) equals $\int_{\lfloor t'\rfloor_{\Delta_{n}}}^{t'}\int_{\lfloor t\rfloor_{\Delta_{n}}}^{t}\partial_{rr'}^{2}p(x,r;y,r')drdr'$,
which can be upper bounded by a double integral, incurring in all
an additional $\Delta_{n}^{2}$. From this the result is obtained.
\end{proof}

\subsection{Proof of Proposition \ref{prop:upperBound_FT}}
\begin{proof}
(i). Let $\mathbb{T}_{T,n}$ and $\tilde{\mathbb{T}}_{n}$ be as in
the proof of Proposition \ref{prop:upperBound_density}. By independence
of $\xi$ we have 
\begin{align}
 & \norm{\Gamma_{T}(f(\cdot+\xi))-\widehat{\Gamma}_{T,n}(f(\cdot+\xi))}_{L^{2}(\P)}^{2}\nonumber\\
 &\quad =\int_{\R^{d}}\E\left[\left(\Gamma_{T}(f(\cdot+x))-\widehat{\Gamma}_{T,n}(f(\cdot+x))\right)^{2}\right]\mu(x)dx\nonumber \\
 & \quad\lesssim\norm{\mu}_{\infty}\int_{\R^{d}}\E\left[\left(\Gamma_{T}(f(\cdot+x))-\widehat{\Gamma}_{T,n}(f(\cdot+x))\right)^{2}\right]dx\label{eq:boundedDensity}\\
 & \quad=\norm{\mu}_{\infty}(2\pi)^{-2d}\int_{\R^{d}}|\c Ff(u)|^{2}\norm{\Gamma_{T}(e^{i\sc u{\cdot}})-\widehat{\Gamma}_{T,n}(e^{i\sc u{\cdot}})}_{L^{2}(\P)}^{2}du,\label{eq:error_char_fun}
\end{align}
using the Plancherel Theorem in the last line. For $u\in\R^{d}$,
$0\leq t,t'\leq T$, set 
\begin{align*}
E_{t,t'}^{u} & =\E[(e^{i\sc u{X_{t}}}-e^{i\sc u{X_{\lfloor t\rfloor_{\Delta_{n}}}}})(e^{i\sc{-u}{X_{t'}}}-e^{i\sc{-u}{X_{\lfloor t'\rfloor_{\Delta_{n}}}}})]\\
 & =\varphi(u,t;-u,t')-\varphi(u,t;-u,\lfloor t'\rfloor_{\Delta_{n}})\\
 & \qquad -\varphi(u,\lfloor t\rfloor_{\Delta_{n}};-u,t')+\varphi(u,\lfloor t\rfloor_{\Delta_{n}};-u,\lfloor t'\rfloor_{\Delta_{n}}).
\end{align*}
Write $\norm{\Gamma_{T}(e^{i\sc u{\cdot}})-\widehat{\Gamma}_{T,n}(e^{i\sc u{\cdot}})}_{L^{2}(\P)}^{2}$
as $\int_{[0,T]^{2}}E_{t,t'}^{u}d(t,t')=A_{1}^{u}+2A_{2}^{u}+2A_{3}^{u}$,
with $A_{i}^{u}$, $i=1,2,3$ as in (\ref{eq:A_1_2_3}) above, but
with $E_{t,t'}$ replaced by $E_{t,t'}^{u}$. For the result it is
enough to show, with
\begin{align*}
w_{T,n}^{(u)} & =\Delta_{n}\int_{\Delta_{n}}^{T}|g(u,t)|dt+T\int_{0}^{\Delta_{n}}|g(u,t)|dt,\\
v_{T,n}^{(u)} & =\Delta_{n}\int_{\mathbb{T}_{T,n}}|\partial_{t'}\varphi(u,t;-u,t')-\partial_{t'}\varphi(u,\lfloor t\rfloor_{\Delta_{n}};-u,t')|d(t,t'),
\end{align*}
that
\begin{align}
|A_{1}^{u}|+|A_{2}^{u}| & \lesssim w_{T,n}^{(u)},\label{eq:A_1_2_charFun}\\
|A_{3}^{u}| & \lesssim v_{T,n}^{(u)}.\label{eq:A_3_charFun}
\end{align}
For $t=t'$, $\overline{E_{t,t}^{u}}=E_{t,t}^{u}$ and so $E_{t,t}^{u}=2-2\text{Re}(\varphi(u,\lfloor t\rfloor_{\Delta_{n}};-u,t))=g(u,t)$.
(\ref{eq:A_1_2_charFun}) follows immediately as for $A_{1}$, $A_{2}$
in the proof of Proposition \ref{prop:upperBound_density}, again
with $E_{t,t'}^{u}$ instead of $E_{t,t'}$. On the other hand, differentiability
of the characteristic functions for $(t,t')\in\mathbb{\tilde{T}}_{n}$
shows
\begin{equation}
E_{t,t'}^{u}=\int_{\lfloor t'\rfloor_{\Delta_{n}}}^{t'}\{\partial_{r}\varphi(u,t;-u,r)-\partial_{r}\varphi(u,\lfloor t\rfloor_{\Delta_{n}};-u,r)\}dr.
\end{equation}
Arguing as after (\ref{eq:ganychenko_trick}) above yields (\ref{eq:A_3_charFun})
and thus the result.

(ii). It is enough to prove (\ref{eq:A_3_charFun}) with $\tilde{v}_{T,n}^{(u)}=\Delta_{n}^{2}\int_{\mathbb{T}_{T,n}}|\partial_{tt'}^{2}\varphi(u,t;-u,t')|d(t,t')$
instead of $v_{T,n}^{(u)}$. As in the proof of Proposition \ref{prop:upperBound_density},
for this it suffices to note by differentiability of $r\mapsto\partial_{r}\varphi(u,r;-u,r')$
that $E_{t,t'}^{u}=\int_{\lfloor t'\rfloor_{\Delta_{n}}}^{t'}\int_{\lfloor t\rfloor_{\Delta_{n}}}^{t}\partial_{rr'}^{2}\varphi(u,r;-u,r')drdr'$.
\end{proof}

\subsection{Proofs of Section \ref{subsec:Fractional-Brownian-motion}}

Observe first the following elementary lemma, which will be used frequently.
\begin{lem}
\label{lem:summations} We have for $\alpha,\beta,\gamma\in\R$: 
\begin{align*}
\int_{\Delta_{n}}^{T}\int_{t+2\Delta_{n}}^{T}\frac{1}{|t'-t|^{\alpha}t^{\beta}(t')^{\gamma}}dt'dt & \lesssim T^{2-\alpha-\beta-\gamma}(\I_{\{\alpha+\gamma=1\}}\log n+\I_{\left\{ \alpha+\gamma\neq1\right\} }\max(1,n^{\alpha+\gamma-1}))\\
 & \quad\cdot(\I_{\{\beta=1\}}\log n+\I_{\left\{ \beta\neq1\right\} }\max(1,n^{\beta-1})).
\end{align*}
\end{lem}

\begin{proof}
By the changes of variables $t=Tr$, $t'=Tr'$ we have
\begin{align*}
\int_{\Delta_{n}}^{T}\int_{t+2\Delta_{n}}^{T}\frac{1}{|t'-t|^{\alpha}t^{\beta}(t')^{\gamma}}dt'dt & =T^{2-\alpha-\beta-\gamma}\int_{1/n}^{1-2/n}\frac{1}{r^{\beta}}\int_{r+2/n}^{1}\frac{1}{|r'-r|^{\alpha}(r')^{\gamma}}dr'dr.
\end{align*}
For fixed $r$ the $dr'$-integral equals
\begin{align*}
\int_{2/n}^{1-r}\frac{1}{(r')^{\alpha}(r'+r)^{\gamma}}dr' & \leq\int_{2/n}^{1}\frac{1}{(r')^{\alpha+\gamma}}dr'\lesssim\I_{\{\alpha+\gamma=1\}}\log n+\I_{\left\{ \alpha+\gamma\neq1\right\} }\max(1,n^{\alpha+\gamma-1}),
\end{align*}
and the same upper bound applies to the $dr$-integral with $\beta$
instead of $\alpha+\gamma$.
\end{proof}

If $X$ is a fractional Brownian motion, then the distribution of
$(X_{t},X_{t'})$, $t\neq t'$, is Gaussian and its characteristic
function is $\varphi(u,t;v,t')=e^{-\frac{1}{2}\Phi_{t,t'}(u,v)}$
with
\[
\Phi_{t,t'}(u,v)=|u|{}^{2}t^{2H}+|v|^{2}(t')^{2H}+2\sc uvc(t,t'),
\]
where we recall the covariance function $c(t,t')$ from (\ref{eq:fracBMcovFun}).
A simple computation shows 
\begin{align}
\partial_{t'}\varphi(u,t;v,t') & =-\frac{1}{2}\partial_{t'}\Phi_{t,t'}(u,v)\varphi(u,t;v,t'),\nonumber \\
\partial_{tt'}^{2}\varphi(u,t;v,t') & =\bigg(-\frac{1}{2}\partial_{tt'}^{2}\Phi_{t,t'}(u,v)+\frac{1}{4}\partial_{t}\Phi_{t,t'}(u,v)\partial_{t'}\Phi_{t,t'}(u,v)\bigg)\varphi(u,t;v,t'),\nonumber \\
\text{with \,\,}\partial_{t'}\Phi_{t,t'}(u,v) & =2H(\sc{u+v}v(t')^{2H-1}-\sc uv|t'-t|^{2H-1}),\nonumber \\
\partial_{t}\Phi_{t,t'}(u,v) & =2H(\sc u{u+v}t^{2H-1}-\sc uv|t'-t|^{2H-1}),\nonumber \\
\partial_{tt'}^{2}\Phi_{t,t'}(u,v) & =2H(2H-1)\sc uv|t'-t|^{2H-2}.\label{eq:Phi}
\end{align}
Self-similarity of a fractional Brownian motion also implies $(X_{t})_{0\leq t\leq T}\overset{d}{\sim}(T^{H}X_{t/T})_{0\leq t\leq T}$,
and therefore 
\begin{align}
\varphi(u,t;v,t') & =\varphi\left(T^{H}u,t/T;T^{H}v,t'/T\right),\label{eq:time_transform}\\
\partial_{tt'}^{2}\Phi_{t,t'}(u,v) & =T^{-2}\partial_{tt'}^{2}\Phi_{t/T,t'/T}\left(T^{H}u,T^{H}v\right),\label{eq:time_transform2}
\end{align}
and similarly for $\partial_{t'}\Phi_{t,t'}(u,v)$, $\partial_{t}\Phi_{t,t'}(u,v)$.
In order to upper bound the characteristic function, we use that a
fractional Brownian motion is \emph{locally nondeterministic}, cf.
\citep{Xiao2006,Berman:1973ke}, that is, for $r'>r$
\begin{align*}
\Phi_{r,r'}(u,v) & =\text{Var}(\sc v{X_{r'}-X_{r}}+\sc{u+v}{X_{r}})\\
 & \geq C(v^{2}|r'-r|^{2H}+(u+v)^{2}r^{2H}).
\end{align*}
In particular, 
\begin{equation}
\varphi(u,r;v,r')\leq e^{-Cv^{2}|r'-r|^{2H}-C(u+v)^{2}r^{2H}}.\label{eq:localNonDet}
\end{equation}
With this let us prove the two theorems and Corollary \ref{cor:localTime}.
\begin{proof}[Proof of Theorem \ref{thm:fractionalBM}]
 It is easily checked by a look at (\ref{eq:Phi}) that the assumptions
of Proposition \ref{prop:upperBound_FT}(i,ii) are satisfied with
$X_{0}=\xi$ $X-X_{0}$ independent of $X_{0}$. For $T=1$ and $v=-u$
we have in (\ref{eq:Phi})
\begin{align*}
\Phi_{t,t'}(u,-u) & =|u|^{2}|t'-t|^{2H},\\
\partial_{t}\Phi_{t,t'}(u,-u) & =\partial_{t'}\Phi_{t,t'}(u,-u)=2H|u|^{2}|t'-t|^{2H-1},\\
\partial_{tt'}^{2}\Phi_{t,t'}(u,-u) & =-2H(2H-1)|u|^{2}|t'-t|^{2H-2},\\
\partial_{tt'}^{2}\varphi(u,t;-u,t') & =\bigg(H(2H-1)|u|^{2}|t'-t|^{2H-2}+H^{2}|u|^{4}|t'-t|^{4H-2}\bigg)e^{-\frac{1}{2}|u|^{2}|t'-t|^{2H}}.
\end{align*}
Consequently, with $g(u,t)=2(1-e^{-\frac{1}{2}|u|^{2}|t-\lfloor t\rfloor_{\Delta_{n}}|^{2H}})$,
\begin{align*}
w_{1,n}^{(u)} & =\Delta_{n}\int_{\Delta_{n}}^{T}|g(u,t)|dt+T\int_{0}^{\Delta_{n}}|g(u,t)|dt\\
 & \lesssim\Delta_{n}\int_{\Delta_{n}}^{1}\left(|u|^{2}|t-\lfloor t\rfloor_{\Delta_{n}}|^{2H}\right)^{s}dt+\int_{0}^{\Delta_{n}}\left(|u|^{2}|t-\lfloor t\rfloor_{\Delta_{n}}|^{2H}\right)^{s}dt\leq\Delta_{n}^{1+2sH}|u|^{2s}.
\end{align*}
On the other hand, by (\ref{eq:Phi}) and $\sup_{x\geq0}(x{}^{2}e^{-x})\lesssim1$,
we find 
\begin{align}
\tilde{v}_{T,n}^{(u)} & =\Delta_{n}^{2}\int_{\mathbb{T}_{T,n}}|\partial_{tt'}^{2}\varphi(u,t;-u,t')|d(t,t')\nonumber \\
 & \lesssim\Delta_{n}^{2}\int_{\mathbb{T}_{1,n}}\bigg(|u|^{2}|t'-t|^{2H-2}+|u|^{4}|t'-t|^{4H-2}\bigg)e^{-\frac{1}{2}|u|^{2}|t'-t|^{2H}}d(t,t')\nonumber \\
 & \lesssim\Delta_{n}^{2}|u|^{2}\int_{\mathbb{T}_{1,n}}|t'-t|^{2H-2}e^{-\frac{1}{2}|u|^{2}|t'-t|^{2H}}d(t,t')\nonumber \\
 & \lesssim\Delta_{n}^{2}|u|^{2}\int_{\Delta_{n}}^{1}\int_{t+2\Delta_{n}}^{1}(t'-t){}^{2H-2}e^{-\frac{1}{2}|u|^{2}(t'-t){}^{2H}}dt'dt\nonumber \\
 & \le\Delta_{n}^{2}|u|^{2}\int_{\Delta_{n}}^{1}(t')^{2H-2}e^{-\frac{1}{2}|u|^{2}(t')^{2H}}dt',\label{eq:fracBM_1}
\end{align}
concluding by the definition of the set $\mathbb{T}_{1,n}$ in (\ref{eq:T_n})
for the last two inequalities. If $2sH<1$, then using $\sup_{x\geq0}(x{}^{\gamma}e^{-x})\lesssim1$
for $\gamma=2sH$ gives
\begin{align*}
\tilde{v}_{T,n}^{(u)} & \lesssim\Delta_{n}^{2}|u|^{2s}\int_{\Delta_{n}}^{1}(t')^{2Hs-2}dt'\lesssim\Delta_{n}^{1+2sH}|u|^{2s},
\end{align*}
while if $2sH=1$, then using in (\ref{eq:fracBM_1}) the change of
variables $|u|^{2}(t')^{1/s}\mapsto z$ shows 
\begin{align*}
\tilde{v}_{T,n}^{(u)} & \lesssim|u|^{2s}\Delta_{n}^{2}s\int_{0}^{|u|^{2}}z^{-s}e^{-\frac{1}{2}z}dz\lesssim\Delta_{n}^{2}|u|{}^{2s}=\Delta_{n}^{1+2sH}|u|{}^{2s},
\end{align*}
implying (\ref{eq:upperBoundH}) for all $0\leq s\leq\min(1,1/(2H))$.
This proves $w_{1,n}^{(u)}+\tilde{v}_{1,n}^{(u)}\lesssim\Delta_{n}^{1+2sH}|u|^{2s}$
and therefore the result for $T=1$. For general $T$ we only note
by the self-similarity identities in (\ref{eq:time_transform}), (\ref{eq:time_transform2})
that $w_{T,n}^{(u)}+\tilde{v}_{T,n}^{(u)}=T^{2}(w_{1,n}^{(T^{H}u)}+\tilde{v}_{1,n}^{(T^{H}u)})$.
\end{proof}
\begin{proof}[Proof of Theorem \ref{thm:fracBM_d_1}]
 Since $X_{0}$ may not have a density, the error of the Riemann
estimator on $[0,2\Delta_{n}]$ can generally not be controlled by
the $L^{2}$-norm of $f$, but for bounded $f$ we have
\[
\norm{\Gamma_{T}(f)-\widehat{\Gamma}_{T,n}(f)}_{L^{2}(\P)}\leq4\norm f_{\infty}\Delta_{n}+R_{n}(f),
\]
with $R_{n}(f)=\norm{\int_{2\Delta_{n}}^{T}(f(X_{t})-f(X_{\lfloor t\rfloor_{\Delta_{n}}}))dt}_{L^{2}(\P)}$.
In order to prove 
\begin{equation}
R_{n}(f)\lesssim\norm f_{H^{s}}^{1/2}\max(\norm f_{H^{s}}^{1/2},\norm f_{L^{2}}^{1/2})\max(1,T^{1/4})T^{1/4-H/2}\Delta_{n}^{1/2+sH},\label{eq:fracBM_d_1_upperBound}
\end{equation}
note that smooth and compactly supported functions are dense in $H^{s}(\R^{d})$.
We can therefore always find a sequence of such functions $f^{(\epsilon)}$
with $\norm{f^{(\epsilon)}-f}_{L^{2}},\norm{f^{(\epsilon)}-f}_{H^{s}}\rightarrow0$
as $\epsilon\rightarrow0$. Moreover, as the distribution of each
$X_{t}$, $t>0$, has a Lebesgue density, this means that also the
error terms converge
\[
\int_{2\Delta_{n}}^{T}(f^{(\epsilon)}(X_{t})-f^{(\epsilon)}(X_{\lfloor t\rfloor_{\Delta_{n}}}))dt\rightarrow\int_{2\Delta_{n}}^{T}(f(X_{t})-f(X_{\lfloor t\rfloor_{\Delta_{n}}}))dt
\]
in $L^{2}(\P)$ and so it is enough to show (\ref{eq:fracBM_d_1_upperBound})
for a smooth and compactly supported $f$. For such $f$ we have by
Fourier inversion $f=(2\pi)^{-1}\int_{\R}\c Ff(u)e^{-iu\cdot}du$
and thus
\begin{align*}
R_{n}(f)^{2} & =(2\pi)^{-2}\int_{\R^{2}}\c Ff(u)\c Ff(v)\int_{2\Delta_{n}}^{T}\int_{2\Delta_{n}}^{T}E_{t,t'}^{u,v}dtdt'd(u,v),
\end{align*}
where for $u,v\in\R$ and (cf. $E_{t,t'}^{u}$ in the proof of Proposition
\ref{prop:upperBound_FT})
\begin{align*}
E_{t,t'}^{u,v} & =\varphi(u,t;v,t')-\varphi(u,t;v,\lfloor t'\rfloor_{\Delta_{n}})-\varphi(u,\lfloor t\rfloor_{\Delta_{n}};v,t')+\varphi(u,\lfloor t\rfloor_{\Delta_{n}};v,\lfloor t'\rfloor_{\Delta_{n}}),
\end{align*}
By symmetry in $u,v$ and $t,t'$ we decompose the previous display
similar to (\ref{eq:A_1_2_3}) into 
\begin{align}
 & (2\pi)^{-2}\int_{\R^{2}}\c Ff(u)\c Ff(v)\int_{[2\Delta_{n},T]^{2}}\I_{\{|t-t'|\leq3\Delta_{n}\}}E_{t,t'}^{u,v}d(t,t')d(u,v)\nonumber \\
 & +(2\pi)^{-2}\int_{\R^{2}}\c Ff(u)\c Ff(v)\int_{2\Delta_{n}}^{T}\int_{t+4\Delta_{n}}^{T}2E_{t,t'}^{u,v}dt'dtd(u,v).\label{eq:fracBM_d_2}
\end{align}
If $t\geq2\Delta_{n}$ and $t'>t+4\Delta_{n}$, then by (\ref{eq:Phi})
\begin{align}
 & \left|\int_{2\Delta_{n}}^{T}\int_{t+4\Delta_{n}}^{T}E_{t,t'}^{u,v}dt'dt\right|=\left|\int_{2\Delta_{n}}^{T}\int_{t+4\Delta_{n}}^{T}\int_{\lfloor t\rfloor_{\Delta_{n}}}^{t}\int_{\lfloor t'\rfloor_{\Delta_{n}}}^{t'}\partial_{rr'}^{2}\varphi(u,r;v,r')dr'drdt'dt\right|\nonumber\\
 & \quad\leq\Delta_{n}^{2}\int_{2\Delta_{n}}^{T}\int_{t+4\Delta_{n}}^{T}\left|\partial_{rr'}^{2}\varphi(u,r;v,r')\right|dr'dr\nonumber\\
 & \quad\lesssim\Delta_{n}^{1+2sH}|uv|^{s}\left(U_{T,n}^{(1)}(u,v)+\sum_{j=1}^{4}U_{T,n}^{(2j)}(u,v)^{1/2}U_{T,n}^{(3j)}(u,v)^{1/2}\right),\label{eq:fracBM_d_1}
\end{align}
with
\begin{align*}
U_{T,n}^{(1)}(u,v)& =|uv|^{-s}\Delta_{n}^{1-2sH}\int_{2\Delta_{n}}^{T}\int_{t+4\Delta_{n}}^{T}|\partial_{rr'}^{2}\Phi_{r,r'}(u,v)|\varphi(u,r;v,r')dr'dr,\nonumber \\
U_{T,n}^{(2j)}(u,v)&=|u|^{-2s}\Delta_{n}^{1-2sH}\int_{2\Delta_{n}}^{T}\int_{t+4\Delta_{n}}^{T}g_{r,r'}^{(2,j)}(u,v)\varphi(u,r;v,r')dr'dr,\nonumber \\
 U_{T,n}^{(3j)}(u,v)&=|v|^{-2s}\Delta_{n}^{1-2sH}\int_{2\Delta_{n}}^{T}\int_{t+4\Delta_{n}}^{T}g_{r,r'}^{(3,j)}(u,v)\varphi(u,r;v,r')dr'dr,
\end{align*}
and with 
\begin{align*}
g_{r,r'}^{(2,1)}(u,v) & =|u|^{2s}|v|^{2-2s}|u+v|^{2}(r')^{4H-2}r^{H-1},\quad g_{r,r'}^{(3,1)}(u,v)=|u|^{2-2s}|v|^{2s}|u+v|^{2}r^{3H-1},\\
g_{r,r'}^{(2,2)}(u,v) & =|u|^{2}|v||r'-r|^{4H-2}r^{-H},\quad g_{r,r'}^{(3,2)}(u,v)=|v+u|^{2}|v|^{3}(r')^{4H-2}r^{H},\\
g_{r,r'}^{(2,3)}(u,v) & =|u|^{3}|u+v||r'-r|^{4H-2},\quad g_{r,r'}^{(3,3)}(u,v)=|u||v|^{2}|u+v|r^{4H-2},\\
g_{r,r'}^{(2,4)}(u,v) & =|u|^{4}|r'-r|^{4H-2},\quad g_{r,r'}^{(3,4)}(u,v)=|v|^{4}|r'-r|^{4H-2}.
\end{align*}
On the other hand, when $t,t'\geq2\Delta_{n}$, $|t-t'|\leq3\Delta_{n}$
and $t\leq\lfloor t'\rfloor_{\Delta_{n}}$, noting that the characteristic
functions $\varphi(u,t;v,t')$ are bounded by $1$ and using (\ref{eq:localNonDet})
with $t\leq\lfloor t'\rfloor_{\Delta_{n}}\leq r'$
\begin{align*}
 & \left|\varphi(u,t;v,t')-\varphi(u,t;v,\lfloor t'\rfloor_{\Delta_{n}})\right|\lesssim\left(\int_{\lfloor t'\rfloor_{\Delta_{n}}}^{t'}\left|\partial_{t'}\varphi(u,t;v,r')\right|dr'\right)^{s}\\
 & \lesssim\left(\int_{\lfloor t'\rfloor_{\Delta_{n}}}^{t'}\left|u+v\right||v|(r')^{2H-1}+|uv||r'-t|^{2H-1})dr'\right)^{s}e^{-C(u+v)^{2}t^{2H}}\\
 & \lesssim\Delta_{n}^{2sH}\left(\left|u+v\right|^{s}|v|^{s}+|uv|^{s}\right)e^{-C(u+v)^{2}\min(t,\lfloor t'\rfloor_{\Delta_{n}})^{2H}},
\end{align*}
The same upper bound holds up to a constant when $t\geq t'$, and
when $\lfloor t'\rfloor_{\Delta_{n}}<t<t'$ by arguing as above on
each of the summands in
\begin{align*}
 & \left|\varphi(u,t;v,t')-\varphi(u,t;v,t)\right|+\left|\varphi(u,t;v,t)-\varphi(u,t;v,\lfloor t'\rfloor_{\Delta_{n}})\right|,
\end{align*}
and consequently,
\begin{align*}
 & \left|\int_{[2\Delta_{n},T]^{2}}\I_{\{|t-t'|\leq3\Delta_{n}\}}E_{t,t'}^{u,v}d(t,t')\right|\\
 & \lesssim\Delta_{n}^{2sH}\int_{[2\Delta_{n},T]}\int_{t}^{t+3\Delta_{n}}\left(\left|u+v\right|^{s}|v|^{s}+|uv|^{s}\right)e^{-C(u+v)^{2}(t-\Delta_{n})^{2H}}dt'dt\\
 & =\Delta_{n}^{1+2sH}\left(|v|^{s}U_{T,n}^{(4)}(u+v)+|uv|^{s}U_{T,n}^{(5)}(u+v)\right),\\
\text{with} & \quad U_{T,n}^{(4)}(u)=\int_{[\Delta_{n},T]}t^{-sH}e^{-Cu^{2}t{}^{2H}}dt,\quad U_{T,n}^{(5)}(u)=\int_{[\Delta_{n},T]}e^{-Cu^{2}t{}^{2H}}dt.
\end{align*}
Combining this with (\ref{eq:fracBM_d_1}), (\ref{eq:fracBM_d_2})
we conclude by the Cauchy-Schwarz inequality
\begin{align*}
 & \Delta_{n}^{-1-2sH}R_{n}(f)^{2}\lesssim\int_{\R^{2}}\left|\c Ff(u)\right|\left|\c Ff(v)\right|\big(|v|^{s}U_{T,n}^{(4)}(u+v)+|uv|^{s}U_{T,n}^{(5)}(u+v)\big)d\left(u,v\right)\\
 & \quad+\int_{\R^{2}}\left|\c Ff(u)\right|\left|\c Ff(v)\right||uv|^{s}\big(U_{T,n}^{(1)}(u,v)+\sum_{j=1}^{4}U_{T,n}^{(2j)}(u,v)^{1/2}U_{T,n}^{(3j)}(u,v)^{1/2}\big)d\left(u,v\right)\\
 & \lesssim\norm f_{H^{s}}\norm f_{L^{2}}\norm{U_{T,n}^{(4)}}_{L^{1}}+\norm f_{H^{s}}^{2}\norm{U_{T,n}^{(5)}}_{L^{1}}\\
 & \quad+\norm f_{H^{s}}^{2}\left(\sup_{u\in\R}\int_{\R}U_{T,n}^{(1)}(u,v)dv\right)^{1/2}\left(\sup_{v\in\R}\int_{\R}U_{T,n}^{(1)}(u,v)du\right)^{1/2}\\
 & \quad+\norm f_{H^{s}}^{2}\sum_{j=1}^{4}\left(\sup_{v\in\R}\int_{\R}U_{T,n}^{(3j)}(u,v)du\right)^{1/2}\left(\sup_{u\in\R}\int_{\R}U_{T,n}^{(2j)}(u,v)dv\right)^{1/2}.
\end{align*}
In order to obtain (\ref{eq:fracBM_d_1_upperBound}) we are therefore
left with showing 
\begin{align}
\norm{U_{T,n}^{(4)}}_{L^{1}},\norm{U_{T,n}^{(5)}}_{L^{1}} & \lesssim\max(1,T^{1/2})T^{1/2-H},\nonumber \\
\sup_{u\in\R}\int_{\R}U_{T,n}^{(i)}(u,v)dv,\sup_{v\in\R}\int_{\R}U_{T,n}^{(i)}(u,v)du & \lesssim T^{1-H},\quad i=1,2,3.\label{eq:fracBM_5}
\end{align}
The first line follows from integrating with respect to $u$. With
respect to the second one, using the self-similarity relation (\ref{eq:time_transform}),
we have
\[
\sup_{u\in\R}\int_{\R}U_{T,n}^{(i)}(u,v)dv=T^{1-H}\sup_{u\in\R}\int_{\R}U_{1,n}^{(i)}(u,v)dv,
\]
and similarly when integrating with respect to $u$. It is therefore
enough to consider only $T=1$. Observe first from (\ref{eq:localNonDet})
and from $\sup_{x\geq0}x^{\beta}e^{-|x|}\lesssim1$ for $\beta\geq0$
and $u,v\in\R^{d}$, using also $|u|\leq|u+v|+|v|$, the inequalities
\begin{align*}
|u+v|^{\beta}\varphi(u,r;v,r') & \lesssim r^{-\beta H}e^{-Cv^{2}|r'-r|^{2H}-C(u+v)^{2}r^{2H}},\\
|u|^{\beta}|v|^{\beta'}\varphi(u,r;v,r') & \lesssim(r^{-\beta H}|r'-r|^{-\beta'H}+|r'-r|^{-(\beta+\beta')H})e^{-Cv^{2}|r'-r|^{2H}-C(u+v)^{2}r^{2H}}.
\end{align*}
Together with $\int_{\R}e^{-C(u+v)^{2}r^{2H}}du$, $\int_{\R}e^{-C(u+v)^{2}r^{2H}}dv$,
$\int_{\R}e^{-Cv^{2}r^{2H}}dv\lesssim r^{-H}$, this implies
\begin{align*}
	& \int_{\R}|uv|^{-s}|\partial_{rr'}^{2}\Phi_{r,r'}(u,v)|\varphi(u,r;v,r')dv\\
 	& \quad\lesssim r^{sH-2H}|r'-r|^{sH+H-2}+r^{-H}|r'-r|^{2sH-2},\\
\end{align*}
which also holds when integrating with respect
to $du$ instead of $dv$, as well as
\begin{align*}
 & \int_{\R}|u|^{-2s}g_{r,r'}^{(2,1)}(u,v)\varphi(u,r;v,r')dv\lesssim|r'-r|^{2sH-3H}(r')^{4H-2}r^{-1-H},\\
 & \int_{\R}|v|^{-2s}g_{r,r'}^{(3,1)}(u,v)\varphi(u,r;v,r')du\lesssim r^{2sH-2H-1}+|r'-r|^{2sH-2H}r^{-1},\\
 & \int_{\R}|u|^{-2s}g_{r,r'}^{(2,2)}(u,v)\varphi(u,r;v,r')dv\lesssim r^{2sH-3H}|r'-r|^{2H-2}+|r'-r|^{2sH-2}r^{-H},\\
 & \int_{\R}|v|^{-2s}g_{r,r'}^{(3,2)}(u,v)\varphi(u,r;v,r')du\lesssim|r'-r|^{2sH-3H}(r')^{4H-2}r^{-2H},\\
 & \int_{\R}|u|^{-2s}g_{r,r'}^{(2,3)}(u,v)\varphi(u,r;v,r')dv\lesssim r^{2sH-4H}|r'-r|^{3H-2}+|r'-r|^{2sH-2}r^{-H}\\
 & \int_{\R}|v|^{-2s}g_{r,r'}^{(3,3)}(u,v)\varphi(u,r;v,r')du\lesssim r^{H-2}|r'-r|^{2sH-2}+|r'-r|^{2sH-3H}r^{2H-2}\\
 & \int_{\R}|u|^{-2s}g_{r,r'}^{(2,4)}(u,v)\varphi(u,r;v,r')dv\lesssim r^{2sH-4H}|r'-r|^{3H-2}+|r'-r|^{2sH-2}r^{-H},\\
 & \int_{\R}|v|^{-2s}g_{r,r'}^{(3,4)}(u,v)\varphi(u,r;v,r')du\lesssim|r'-r|^{2sH-2}r^{-H}.
\end{align*}
Plugging these estimates into the definitions
of the $U_{1,n}^{(i)}(u,v)$ with $sH\leq H$, $H\leq1/2$ and using
Lemma \ref{lem:summations} yields (\ref{eq:fracBM_5}), and thus
finishes the proof.
\end{proof}
\begin{proof}[Proof of Corollary \ref{cor:occupationTime}]
We use an approximation argument similar to Theorem 3.6 of \citep{Altmeyer:2017dza}
for $f=\I_{[a,b]}$. Let first $H>1/2$ and let $a,b$ be finite.
Consider the smooth convolution approximation $f_{\epsilon}=f*\varphi_{\epsilon}$,
where $\varphi_{\epsilon}(x)=\epsilon^{-1}\varphi(\epsilon^{-1}x)$
for $\epsilon>0$ and a smooth $\varphi$ with support in $[-1,1]$,
$\int_{\R}\varphi(x)dx=1$. Then $f_{\epsilon}\in H^{1}(\R^{d})$
and for sufficiently small $\epsilon$, $f-f_{\epsilon}$ is supported
in $\text{\ensuremath{\c A}=}[a-\epsilon,a+\epsilon]\cup[b-\epsilon,b+\epsilon]$.
Consequently, $\norm{f-f_{\epsilon}}_{L^{2}}^{2}\lesssim\epsilon$
and
\begin{align*}
\norm{f_{\epsilon}}_{H^{1}}^{2} & =\epsilon^{-1}\int\left|\frac{f(x)}{\epsilon}\int\varphi'(y)dy-f_{\epsilon}'(x)\right|^{2}dx\\
 & =\epsilon^{-2}\int_{\c A}\left|\int(f(x)-f(x+\epsilon y))\varphi'(y)dy\right|^{2}dx\lesssim\epsilon^{-1},
\end{align*}
and the implied constants are independent of $a,b$. By Theorem \ref{thm:fractionalBM}
with $s=0$ and $s=1$ this means
\begin{align*}
 & \norm{\Gamma_{T}(f)-\widehat{\Gamma}_{T,n}(f)}_{L^{2}(\P)}\leq\norm{\Gamma_{T}(f-f_{\epsilon})-\widehat{\Gamma}_{T,n}(f-f_{\epsilon})}_{L^{2}(\P)}+\norm{\Gamma_{T}(f_{\epsilon})-\widehat{\Gamma}_{T,n}(f_{\epsilon})}_{L^{2}(\P)}\\
 & \quad\lesssim\norm{\mu}_{\infty}^{1/2}\left(\norm{f-f_{\epsilon}}_{L^{2}}T^{1/2}\Delta_{n}^{1/2}+\norm{f_{\epsilon}}_{H^{s}}T^{1/2}\Delta_{n}^{1/2+H}\right)\\
 & \quad\lesssim\norm{\mu}_{\infty}^{1/2}T^{1/2}\Delta_{n}^{1/2}\left(\epsilon+\epsilon^{-1}\Delta_{n}^{H}\right).
\end{align*}
The claim for finite $a,b$ follows with $\epsilon=\Delta_{n}^{H/2}$.
Since the upper bound is independent of $a,b$, letting $a\rightarrow-\infty$
and $b\rightarrow\infty$ yields the claim for all $a,b$ and $H>1/2$.
The proof for $H\leq1/2$ follows in the same way from Theorem \ref{thm:fracBM_d_1},
noting $\norm f_{\infty},\norm{f_{\epsilon}}_{\infty}\lesssim1$,
but this time $\norm f_{L^{2}}$ and thus the implied constants above
depend on $a,b$.
\end{proof}
\begin{proof}[Proof of Corollary \ref{cor:localTime}]
 We have from Corollary \ref{cor:occupationTime}
\begin{align*}
\lVert L_{T}(a)-\widehat{\Gamma}_{T,n}(f_{a,n})\rVert_{L^{2}(\P)} & \lesssim\norm{L_{T}(a)-\Gamma_{T}(f_{a,n})}_{L^{2}(\P)}+\max(1,T^{1/4})T^{1/4-H/2}\Delta_{n}^{1/2-H/2}.
\end{align*}
Self-similarity of $X$ implies that $L_{T}$ has the same distribution
as $T^{1-H}L_{1}(T^{-H}\cdot)$. By the occupation time formula, cf.
\citep{Geman1980}, this means also $\int_{0}^{T}f_{a,n}(X_{t})dt=\int_{\R}f_{a,n}(x)L_{T}(x)dx$
and $T^{1-H}\int_{\R}f_{a,n}(x)L_{1}(T^{-H}x)dx$ have the same distribution.
With $\int_{\R}f_{a,n}(x)dx=1$,
\begin{align*}
 & \norm{L_{T}(a)-\Gamma_{T}(f_{a,n})}_{L^{2}(\P)}\lesssim T^{1-H}\int_{-1}^{1}\norm{L_{1}(T^{-H}a)-L_{1}(T^{-H}(\Delta_{n}^{H}x+a))}_{L^{2}(\P)}dx\\
 & \lesssim T^{1-H}\int_{-1}^{1}(T^{-H}\Delta_{n}^{H}x)^{\gamma}dx\lesssim T^{1-H-\gamma H}\Delta_{n}^{\gamma H},
\end{align*}
uniformly in $a\in\R$ with $\gamma<(1-H)/(2H)$, where the last line
follows from moment bounds for the local time (e.g., \citep{Berman:1969jm}
or Equation 4.18 of \citep{Xiao2006}). Choosing $\gamma$ arbitrarily
close to $(1-H)/(2H)$ gives the result.
\end{proof}

\subsection{Proof of Theorem \ref{thm:density_L2_Hs}}
\begin{proof}
Without loss of generality we can assume $f$ to be smooth and compactly
supported. Indeed, we can always find a sequence of smooth and compactly
supported functions $f^{(\epsilon)}$ with $\norm{f^{(\epsilon)}-f}_{L^{2}},\norm{f^{(\epsilon)}-f}_{H^{s}}\rightarrow0$
as $\epsilon\rightarrow0$. Since the distribution of each $X_{t}$,
$t\geq0$, has a Lebesgue density, this means $\Gamma_{T}(f^{(\epsilon)})\rightarrow\Gamma_{T}(f)$
and $\widehat{\Gamma}_{T,n}(f^{(\epsilon)})\rightarrow\widehat{\Gamma}_{T,n}(f)$
in $L^{2}(\P)$ and so the claim of the theorem transfers from $f^{(\epsilon)}$
to $f$.

The density of $(X_{t},X_{t'})$ for $t<t'$, $x,y\in\R^{d}$, is
$p(t,x;t',y)=\int p(t,x;t',y;x_{0})\mu(x_{0})dx_{0}$. As $\mu$ is
bounded, 
\[
\int_{\R^{d}}q_{t'}(x-x_{0})q_{t'-t}(y-x)\mu(x_{0})dx_{0}\leq\norm{\mu}_{\infty}q_{t'-t}(y-x).
\]
The respective heat kernel bounds from Assumptions \ref{ass:A} and
\ref{ass:B} yield then, using (\ref{eq:heat_bound_1}), (\ref{eq:heat_bound_2}),
(\ref{eq:heat_bound_3}) above,
\begin{align}
p(x,t;y,t') & \lesssim\norm{\mu}_{\infty}q_{t'-t}(y-x),\label{eq:heat_bound_1-1}\\
|\partial_{t'}p(x,t;y,t')| & \lesssim\norm{\mu}_{\infty}\frac{1}{t'-t}q_{t'-t}(y-x),\label{eq:heat_bound_2-1}\\
|\partial_{tt'}^{2}p(x,t;y,t')| & \lesssim\norm{\mu}_{\infty}(\frac{1}{t(t'-t)}+\frac{1}{(t'-t)^{2}})q_{t'-t}(y-x).\label{eq:heat_bound_3-1}
\end{align}
This means
\begin{align*}
w_{T,n}^{(x,y)} & \lesssim\Delta_{n}\int_{\Delta_{n}}^{T}p(x,\lfloor t\rfloor_{\Delta_{n}};y,t)dt+T\int_{0}^{\Delta_{n}}p(x,\lfloor t\rfloor_{\Delta_{n}};y,t)dt\\
 & \lesssim\norm{\mu}_{\infty}\Delta_{n}\int_{\Delta_{n}}^{T}q_{t'-t}(y-x)dt+T\int_{0}^{\Delta_{n}}q_{t'-t}(y-x)dt,\\
v_{T,n}^{(x,y)} & \lesssim\Delta_{n}\int_{\mathbb{T}_{T,n}}|\partial_{t'}p(x,t;y,t')-\partial_{t'}p(x,\lfloor t\rfloor_{\Delta_{n}};y,t')|d(t,t')\\
 & \lesssim\norm{\mu}_{\infty}\Delta_{n}\int_{\mathbb{T}_{T,n}}\frac{1}{t'-t}(q_{t'-t}(y-x)+q_{t'-\lfloor t\rfloor_{\Delta_{n}}}(y-x))d(t,t'),\\
\tilde{v}_{T,n}^{(x,y)} & \lesssim \norm{\mu}_{\infty \Delta_{n}^{2}\int_{\mathbb{T}_{T,n}}}(\frac{1}{t(t'-t)}+\frac{1}{(t'-t)^{2}})q_{t'-t}(y-x)d(t,t').
\end{align*}
Moreover, (\ref{eq:densityFourierUpper}) and Assumption \ref{ass:B}
with $\rho=q_{t'-t}$ show for $0\leq s\leq\alpha/2$
\begin{align}
 & \int_{\R^{2d}}(f(x)-f(y))^{2}q_{t'-t}(y-x)d(x,y)\label{eq:f_H_s}\\
 & \quad\lesssim\norm f_{H^{s}}^{2}\int_{\R^{d}}|y|^{2s}q_{t'-t}(y)dy\lesssim\norm f_{H^{s}}^{2}|t'-t|^{2s/\alpha}.\nonumber 
\end{align}
Proposition \ref{prop:upperBound_density}(i,ii) then yields for $s=0$
\begin{align*}
\norm{\Gamma_{T}(f)-\widehat{\Gamma}_{T,n}(f)}_{L^{2}(\P)}^{2} & \lesssim\int_{\R^{2d}}(f(x)-f(y))^{2}(w_{T,n}^{(x,y)}+v_{T,n}^{(x,y)})d(x,y)\\
 & \lesssim\norm{\mu}_{\infty}^{2}\norm f_{L^{2}}^{2}(\Delta_{n}T+\Delta_{n}\int_{\mathbb{T}_{T,n}}\frac{1}{t'-t}d(t,t')),
\end{align*}
while for $s>0$ we get instead the upper bound
\begin{align*}
 & \norm{\mu}_{\infty}^{2}\norm f_{H^{s}}^{2}(T\Delta_{n}^{1+2s/\alpha}
 +\Delta_{n}^{2}\int_{\mathbb{T}_{T,n}}(\frac{1}{t(t'-t)^{1-2s/\alpha}}+\frac{1}{(t'-t)^{2-2s/\alpha}})d(t,t')).
\end{align*}
The result follows from Lemma \ref{lem:summations}.
\end{proof}

\subsection{Proof of Theorem \ref{thm:density_bdd_Holder}}
\begin{proof}[Proof of Theorem \ref{thm:density_bdd_Holder}]
 Under Assumptions \ref{ass:A} and \ref{ass:B}, respectively, the
required integrabilities of $\partial_{t'}p(x,t;y,t';x_{0})$ and
$\partial_{tt'}^{2}p(x,t;y,t';x_{0})$ follow from continuity of $\partial_{t'}p_{t,t'}(x,y)$
and $\partial_{tt'}^{2}p_{t,t'}(x,y)$ on $\mathbb{T}_{T,n}$. Formally,
we have for $t<t'$ 
\begin{align*}
\partial_{t'}p(x,t;y,t';x_{0}) & =p_{0,t}(x_{0},x)\partial_{t'}p_{t,t'}(x,y),\\
\partial_{tt'}^{2}p(x,t;y,t';x_{0}) & =\partial_{t}p_{0,t}(x_{0},x)\partial_{t'}p_{t,t'}(x,y)+p_{0,t}(x_{0},x)\partial_{tt'}^{2}p_{t,t'}(x,y).
\end{align*}
Proposition \ref{prop:upperBound_density}(i,ii) yields with $\P_{x_{0}}$
and $p(x,t;y,t';x_{0})$ instead of $\P$ and $p(x,t;y,t')$ for bounded
$f$ or $f\in C^{s}(\R^{d})$, respectively, that 
\begin{align*}
\norm{\Gamma_{T}(f)-\widehat{\Gamma}_{T,n}(f)}_{L^{2}(\P_{x_{0}})}^{2} & \lesssim\int_{\R^{2d}}(f(x)-f(y))^{2}(w_{T,n}^{(x,y)}+v_{T,n}^{(x,y)})d(x,y)\\
 & \leq4\norm f_{\infty}^{2}\int_{\R^{2d}}(w_{T,n}^{(x,y)}+v_{T,n}^{(x,y)})d(x,y),\\
\norm{\Gamma_{T}(f)-\widehat{\Gamma}_{T,n}(f)}_{L^{2}(\P_{x_{0}})}^{2} & \lesssim\norm f_{C^{s}}^{2}\int_{\R^{2d}}|x-y|^{2s}(w_{T,n}^{(x,y)}+\tilde{v}_{T,n}^{(x,y)})d(x,y).
\end{align*}
The heat kernel bounds on $p_{t,t'}$ and the formal derivatives of
$p(x,t;y,t';x_{0})$ above show 
\begin{align}
p(x,t;y,t';x_{0}) & \lesssim q_{t}(x-x_{0})q_{t'-t}(y-x),\label{eq:heat_bound_1}\\
|\partial_{t'}p(x,t;y,t';x_{0})| & \lesssim\frac{1}{t'-t}q_{t}(x-x_{0})q_{t'-t}(y-x),\label{eq:heat_bound_2}\\
|\partial_{tt'}^{2}p(x,t;y,t';x_{0})| & \lesssim(\frac{1}{t(t'-t)}+\frac{1}{(t'-t)^{2}})q_{t}(x-x_{0})q_{t'-t}(y-x).\label{eq:heat_bound_3}
\end{align}
Recall that the $q_{t}$ are probability densities. Using (\ref{eq:heat_bound_1})
and (\ref{eq:heat_bound_2}) we have 
\begin{align*}
 & \int_{\R^{2d}}(w_{T,n}^{(x,y)}+v_{T,n}^{(x,y)})d(x,y)\\
 & \quad\lesssim\Delta_{n}\int_{\Delta_{n}}^{T}\int_{\R^{2d}}p(x,\lfloor t\rfloor_{\Delta_{n}};y,t)d(x,y)dt+T\int_{0}^{\Delta_{n}}\int_{\R^{2d}}p(x,\lfloor t\rfloor_{\Delta_{n}};y,t)d(x,y)dt\\
 & \qquad+\Delta_{n}\int_{\mathbb{T}_{T,n}}\int_{\R^{2d}}|\partial_{t'}p(x,t;y,t')-\partial_{t'}p(x,\lfloor t\rfloor_{\Delta_{n}};y,t')|d(x,y)d(t,t')\\
 & \quad\lesssim T\Delta_{n}+\Delta_{n}\int_{\mathbb{T}_{T,n}}\frac{1}{t'-t}d(t,t')\lesssim T\Delta_{n}\log n,
\end{align*}
concluding by Lemma \ref{lem:summations} in the last inequality.
This proves the first part of the result, that is, when $f$ is only
bounded. For the second part set
\begin{align*}
h(t,t') & :=\int_{\R^{2d}}|y-x|^{2s}q_{t}(x-x_{0})q_{t'-t}(y-x)d(x,y).
\end{align*}
Under Assumption \ref{ass:B} we have $h(t,t')\leq\int_{\R^{2d}}|y|^{2s}q_{t'-t}(y)dy\lesssim|t'-t|^{2s/\alpha}$.
Combining this with (\ref{eq:heat_bound_1}), (\ref{eq:heat_bound_2})
and Lemma \ref{lem:summations} yields 
\begin{align*}
 & \int_{\R^{2d}}|y-x|^{2s}(w_{T,n}^{(x,y)}+\tilde{v}_{T,n}^{(x,y)})d(x,y)
 \lesssim\Delta_{n}\int_{\Delta_{n}}^{T}h(\lfloor t\rfloor_{\Delta_{n}},t)dt+T\int_{0}^{\Delta_{n}}h(\lfloor t\rfloor_{\Delta_{n}},t)dt\\
 & \qquad+\Delta_{n}^{2}\int_{\mathbb{T}_{T,n}}\frac{1}{t(t'-t)}+\frac{1}{(t'-t)^{2}})h(t,t')d(t,t')\\
 & \quad\lesssim T\Delta_{n}^{1+2s/\alpha}+\Delta_{n}^{2}\int_{\mathbb{T}_{T,n}}(\frac{1}{t(t'-t)^{1-2s/\alpha}}+\frac{1}{(t'-t)^{2-2s/\alpha}})d(t,t')\\
 & \quad\lesssim T\Delta_{n}^{1+2s/\alpha}+T\Delta_{n}^{1+2s/\alpha}n^{2s/\alpha-1}(\log n+\I_{\{s=0\}}(\log n)^{2}+\I_{\{2s/\alpha=1\}}\log n)\\
 & \quad\lesssim T\Delta_{n}^{1+2s/\alpha}(1+\I_{\{2s/\alpha=1\}}\log n),
\end{align*}
from which we obtain the second part of the result.
\end{proof}

\subsection{Proof of Theorem \ref{thm:additiveProcess}}
\begin{proof}
(i). The characteristic functions $\varphi(\cdot,t;\cdot,t')$ of
$X$ and $X-X_{0}$ evaluated at $(u,-u)$ coincide and so the assumptions
of Proposition \ref{prop:upperBound_FT}(ii) are satisfied with $\xi=X_{0}$.
For $0<t<t'$ we have $\varphi(u,t;-u,t')=e^{\Psi_{t,t'}(-u)}$,
\begin{align*}
\partial_{t'}\Psi_{t,t'}(-u) & =-i\sc u{b_{t'}}-\frac{1}{2}|\sigma_{t'}^{\top}u|^{2}+\int_{\R^{d}}(e^{-i\sc ux}-1+i\sc ux\I_{\left\{ |x|\leq1\right\} })dF_{t'}(x),
\end{align*}
and $\partial_{tt'}^{2}\Psi_{t,t'}(-u)=0$. With $g(u,t)$ from Proposition
\ref{prop:upperBound_FT}(i) and $\gamma=2s/\alpha^{*}\leq1$, Assumption
\ref{assu:II} shows
\begin{align*}
|g(u,t)| & \lesssim|\Psi_{\lfloor t\rfloor_{\Delta_{n}},t}(-u)|^{\gamma}\lesssim\max(1,|u|^{2s})\Delta_{n}^{2s/\alpha^{*}},\\
|\varphi(u,t;-u,t')| & =|e^{\Psi_{t,t'}(-u)}|\lesssim e^{-C|u|^{\alpha}(t'-t)},\\
\partial_{t'}\varphi(u,t;-u,t') & =e^{\Psi_{t,t'}(-u)}\partial_{t'}\Psi_{t,t'}(-u),\\
|\partial_{tt'}^{2}\varphi(u,t;-u,t')| & =|e^{\Psi_{t,t'}(-u)}\partial_{t}\Psi_{t,t'}(-u)\partial_{t'}\Psi_{t,t'}(-u)|\lesssim\max(1,|u|^{2\alpha^{*}})e^{-C|u|^{\alpha}(t'-t)}.
\end{align*}
This yields 
\begin{align*}
w_{T,n}^{(u)} & =\Delta_{n}\int_{\Delta_{n}}^{T}|g(u,t)|dt+T\int_{0}^{\Delta_{n}}|g(u,t)|dt\lesssim\max(1,|u|^{2s})T\Delta_{n}^{1+2s/\alpha^{*}}.
\end{align*}
For $|u|\geq1$ with $\gamma'=2(\alpha^{*}-s)/\alpha$ we have $(|u|^{\alpha}|t'-t|)^{\gamma'}e^{-C|u|^{\alpha}(t'-t)}\lesssim1$
and thus
\begin{align*}
\tilde{v}_{T,n}^{(u)} & =\Delta_{n}^{2}\int_{\mathbb{T}_{T,n}}|\partial_{tt'}^{2}\varphi(u,t;-u,t')|d(t,t')\lesssim|u|{}^{2\alpha^{*}}\Delta_{n}^{2}\int_{\mathbb{T}_{T,n}}e^{-C|u|^{\alpha}(t'-t)}d(t,t')\\
 & \lesssim|u|^{2s}\Delta_{n}^{2}\int_{\mathbb{T}_{T,n}}|t'-t|^{-2(\alpha^{*}-s)/\alpha}d(t,t')\\
 & \lesssim|u|^{2s}\Delta_{n}^{2}T^{2-2(\alpha^{*}-s)/\alpha}(\I_{\{2(\alpha^{*}-s)=\alpha\}}\log n+n^{2(\alpha^{*}-s)/\alpha-1})\\
 & \lesssim|u|^{2s}T\Delta_{n}^{2}(1+\I_{\{2(\alpha^{*}-s)=\alpha\}}\log n),
\end{align*}
using in the last two lines Lemma \ref{lem:summations} and because
always $2(\alpha^{*}-s)\alpha-1\geq0$ for $s\leq\alpha^{*}/2$. By
a different argument for $2(\alpha^{*}-s)=\alpha$ the $\log$-term
can be removed. Indeed, upper bounding $\tilde{v}_{T,n}^{(u)}$ and
integrating over $t'$ in that case yields
\begin{align*}
\tilde{v}_{T,n}^{(u)} & \lesssim|u|^{2\alpha^{*}}\Delta_{n}^{2}\int_{0}^{T}\int_{t}^{T}e^{-C|u|^{\alpha}(t'-t)}d(t,t')\\
 & \lesssim|u|^{2\alpha^{*}-\alpha}\Delta_{n}^{2}\int_{0}^{T}(1-e^{-C|u|^{\alpha}(T-t)})dt\lesssim|u|^{2s}T\Delta_{n}^{2}.
\end{align*}
The same estimates show for $|u|\leq1$ that $\tilde{v}_{T,n}^{(u)}\lesssim T\Delta_{n}^{2}$,
because $2\alpha^{*}-\alpha\geq0$. The result follows from this,
the upper bound on $w_{T,n}^{(u)}$ and (\ref{eq:fourierUpper}) with
$\bar{\rho}=w_{T,n}^{(u)}+v_{T,n}^{(u)}$ such that by Proposition
\ref{prop:upperBound_FT}(ii)
\begin{align*}
\norm{\Gamma_{T}(f)-\widehat{\Gamma}_{T,n}(f)}_{L^{2}(\P)}^{2} & \lesssim\norm{\mu}_{\infty}\int_{\R^{d}}|\c Ff(u)|^{2}(w_{T,n}^{(u)}+v_{T,n}^{(u)})du\\
 & \lesssim\norm{\mu}_{\infty}(\int_{\R^{d}}|\c Ff(u)|^{2}((1+|u|^{2s})T\Delta_{n}^{1+2s/\alpha^{*}}+|u|^{2s}T\Delta_{n}^{2})du\\
 & \lesssim\norm{\mu}_{\infty}\max(\norm f_{L^{2}}^{2},\norm f_{H^{s}}^{2})T\Delta_{n}^{1+2s/\alpha^{*}}.
\end{align*}
(ii). With $|\Psi_{t,t'}(-u)|\lesssim|t'-t|$ we have this time $|g(u,t)|\lesssim\Delta_{n}$
and
\begin{align*}
|\varphi(u,t;-u,t')|,|\partial_{tt'}^{2}\varphi(u,t;-u,t')| & \lesssim|e^{\Psi_{t,t'}(-u)}|.
\end{align*}
Since this is bounded, we immediately find as in (i), $w_{T,n}^{(u)}\lesssim T\Delta_{n}^{2}$,
$\tilde{v}_{T,n}^{(u)}\lesssim T^{2}\Delta_{n}^{2}$. For the supplement
it is enough to note that $|e^{\Psi_{t,t'}(-u)}|\leq e^{-C(t'-t)}$
such that $\tilde{v}_{T,n}^{(u)} \lesssim\Delta_{n}^{2}\int_{\mathbb{T}_{T,n}}e^{-C(t'-t)}d(t,t')\lesssim T\Delta_{n}^{2}$.
\end{proof}

\subsection{Proofs of Section \ref{sec:Lower-bounds}}
\begin{proof}[Proof of Theorem \ref{thm:lowerBound_Hs}]
By definition of $g$ we have $g(1+|u|)^{s/2}\in L^{2}(\R^{d})$
and therefore $f=f^{*}\in H^{s}(\R^{d})$. Assume $T>0$. The sigma
field $\c G_{n}$ is generated by $X_{0}$ and the increments $X_{t_{k}}-X_{t_{k-1}}$
for $k\in\{1,\dots,n\}$. Since they are independent, the Markov property
shows $\E[f(X_{t})|\c G_{n}]=\E[f(X_{t})|X_{t_{k-1}},X_{t_{k}}]$.
This also shows that the random variables $Y_{k}=\int_{t_{k-1}}^{t_{k}}(f(X_{t})-\E[f(X_{t})|\c G_{n}])dt$
are uncorrelated, implying
\[
\norm{\Gamma_{T}(f)-\E\left[\l{\Gamma_{T}(f)}\c G_{n}\right]}_{L^{2}(\P)}^{2}=\sum_{k=1}^{n}\E\left[Y_{k}^{2}\right]=\sum_{k=1}^{n}\E\left[\text{Var}_{k}\left(\int_{t_{k-1}}^{t_{k}}f(X_{t})dt\right)\right],
\]
where $\text{Var}_{k}(Z)$ denotes the variance of a random variable
$Z$, conditional on the sigma algebra $\sigma(X_{t_{k-1}},X_{t_{k}})=\sigma(X_{t_{k-1},}X_{t_{k}}-X_{t_{k-1}})$.
Stationarity and independence of increments yield
\begin{align}
\E\left[\text{Var}_{k}\left(\int_{t_{k-1}}^{t_{k}}f(X_{t})dt\right)\right] & =\Delta_{n}^{2}\E\left[\text{Var}_{0}\left(\int_{0}^{1}f(\Delta_{n}^{1/2}B_{t}+X_{t_{k-1}})dt\right)\right]\nonumber \\
 & =\Delta_{n}^{2}\E\left[\int_{\R^{d}}h_{n}(x)^{2}p_{t_{k-1}}(x)dx\right]\label{eq:lowerBound_4}
\end{align}
for another independent $d$-dimensional Brownian motion $B=(B_{t})_{0\leq t\leq1}$,
with $p_{t}$ denoting the marginal density of $X_{t}$ and with a
random function
\[
h_{n}=\int_{0}^{1}\left(f(\Delta_{n}^{1/2}B_{t}+\cdot)-\E\left[\l{f(\Delta_{n}^{1/2}B_{t}+\cdot)}B_{1}\right]\right)dt.
\]
For $k\geq1$ and $0<T_{0}<t_{1}$ we have in all
\begin{align}
 & \log(\Delta_{n}^{-1/2})^{2}\Delta_{n}^{-(1+s)}\norm{\Gamma_{T}(f)-\E\left[\l{\Gamma_{T}(f)}\c G_{n}\right]}_{L^{2}(\P)}^{2}\gtrsim\E\left[V_{n}\right],\label{eq:lowerBound_2}\\
 & \,\text{with }V_{n}=\int_{\R^{d}}\Delta_{n}^{-d/2}\bar{h}_{n}(\Delta_{n}^{-1/2}x)^{2}\bar{p}_{n}(x)dx,\label{eq:V_n}\\
 & \qquad\,\bar{h}_{n}(x)=\log(\Delta_{n}^{-1/2})\Delta_{n}^{-s/2-d/4}h_{n}(\Delta_{n}^{1/2}x),\quad\bar{p}_{n}(x)=\Delta_{n}\sum_{k=\lceil n/2\rceil}^{n}p_{t_{k-1}}(x).\nonumber 
\end{align}
For $u\in\R^{d}$ set $\kappa(u)=\int_{0}^{1}(e^{-i\sc u{B_{t}}}-\E[e^{-i\sc u{B_{t}}}|B_{1}])dt$.
We will show below
\begin{equation}
0<\E\left[\int_{\R^{d}}|\c F\bar{h}(u)|^{2}du\right]<\infty,\label{eq:bar_h_0}
\end{equation}
with $\bar{h}=\c F^{-1}(|\cdot|^{-s-d/2}\kappa(\cdot))$, which is
therefore $\P$-almost surely well-defined in $L^2(\R^d)$. Recall that $\c Ff=g$
with $g$ from the statement of the Theorem. By properties of the
Fourier transform it follows
\begin{align*}
 & \c F\bar{h}_{n}(u)=\log(\Delta_{n}^{-1/2})\Delta_{n}^{-s/2-d/4}\c Fh_{n}(\Delta_{n}^{-1/2}u)=\frac{\log(\Delta_{n}^{-1/2})\kappa(u)}{(\Delta_{n}^{1/2}+|u|)^{s+d/2}\log(e+\Delta_{n}^{-1/2}|u|)}.
\end{align*}
Writing $\log(e+\Delta_{n}^{-1/2}|u|)=\log(\Delta_{n}^{1/2}e+|u|)+\log(\Delta_{n}^{-1/2})$,
we see that $\log(\Delta_n^{-1/2})/\log(e+\Delta_n^{-1/2}|u|)$ converges to $1$ as $n\rightarrow\infty$ for fixed $u\neq0$, and is for $|u|>1$
also upper bounded by $1$, while for $|u|\leq1$ this follows
from $\log(\Delta_{n}^{-1/2})=\log(\Delta_{n}^{-1/2}|u|)-\log(|u|)\leq\log(e+\Delta_{n}^{-1/2}|u|)$.
In particular, we have $\P$-almost surely for all $u\in\R^{d}$
\begin{align}
 & \c F\bar{h}_{n}(u)\r{}\c F\bar{h}(u),\quad n\rightarrow\infty,\label{eq:bar_h}\\
 & \left|\c F\bar{h}_{n}(u)\right|\leq|\c F\bar{h}(u)|.\label{eq:bar_h_2}
\end{align}
The dominated convergence theorem implies therefore $\P$-almost surely
the $L^{2}$-convergence $\c F\bar{h}_{n}\rightarrow\c F\bar{h}$
or equivalently $\bar{h}_{n}\rightarrow\bar{h}$. Moreover, $\bar{p}_{n}\r{}\int_{T/2}^{T}p_{t}dt$
in $L^{1}(\R^{d})$ and $\sup_{x\in\R^{d}}\bar{p}_{n}(x)$ is bounded.
By a change of variables and with $V_{n}$ from (\ref{eq:V_n}) this
means $\P$-almost surely
\begin{align*}
 & V_{n}=\int_{\R^{d}}\bar{h}_{n}(x)^{2}\bar{p}_{n}(\Delta_{n}^{1/2}x)dx\rightarrow\int_{T/2}^{T}p_{t}(0)dt\int_{\R^{d}}\bar{h}(x)^{2}dx,\\
 & V_{n}\lesssim\int_{\R^{d}}\bar{h}_{n}(x)^{2}dx\lesssim\int_{\R^{d}}|\bar{h}(u)|^{2}du,
\end{align*}
concluding by the Plancherel theorem and (\ref{eq:bar_h_2}) in the
last line. The result follows from using the dominated convergence
theorem in (\ref{eq:lowerBound_2}) and from (\ref{eq:bar_h_0}),
$\int_{T/2}^{T}p_{t}(0)\gtrsim T^{1-d/2}$.

We are left with proving (\ref{eq:bar_h_0}). The lower bound is clearly
true, since $\kappa$ does not vanish identically. For the upper bound
it is enough by Fubinis Theorem to show
\begin{equation}
\E\left[|\kappa(u)|^{2}\right]\lesssim\min(|u|^{2},1),\quad u\in\R^{d}.\label{eq:lowerBound_3}
\end{equation}
Since the complex exponentials $e^{-i\sc u{B_{t}}}$ are bounded in
absolute value, we clearly have $\E[|\kappa(u)|^{2}]\lesssim1$ for
all $u\in\R^{d}$. On the other hand, if $|u|\leq1$, then 
\begin{align*}
& \E\left[|\kappa(u)|^{2}\right] =\E\left[\left|\int_{0}^{1}(e^{-i\sc u{B_{t}}}-1)dt-\E[\int_{0}^{1}(e^{-i\sc u{B_{t}}}-1)dt|B_{1}])\right|^{2}\right]\\
 & \quad \leq\E\left[\left|\int_{0}^{1}\left(e^{-i\sc u{B_{t}}}-1\right)dt\right|\right]^{2}\leq\int_{0}^{1}\E\left[\left|e^{-i\sc u{B_{t}}}-1\right|^{2}\right]dt\leq|u|^{2}\int_{0}^{1}|B_{t}|^{2}dt\lesssim|u|^{2},
\end{align*}
implying the wanted upper bound in (\ref{eq:lowerBound_3}).
\end{proof}
\begin{proof}[Proof of Corollary \ref{cor:lowerBound}]
The proof follows along the lines of Theorem \ref{thm:lowerBound_Hs}.
Since $X_{0}$ is independent of $(X_{t}-X_{0})_{t\geq0}$, we can
argue in (\ref{eq:lowerBound_4}) with $X_{t}-X_{0}$ instead of $X_{t}$,
such that
\begin{align*}
 & \E\left[\text{Var}_{k}\left(\int_{t_{k-1}}^{t_{k}}f(X_{t})dt\right)\right]=\E\left[\text{Var}_{k}\left(\int_{t_{k-1}}^{t_{k}}f(X_{t}-X_{0}+X_{0})dt\right)\right]\\
 & \quad=\Delta_{n}^{2}\int_{\R^{d}}\E\left[\int_{\R^{d}}h_{n}(x+y)^{2}p_{t_{k-1}}(x)dx\right]\mu(y)dy\\
 & \quad=\Delta_{n}^{2}\E\left[\int_{\R^{d}}h_{n}(x)^{2}\left(\int_{\R^{d}}p_{t_{k-1}}(x-y)\mu(y)dy\right)dx\right].
\end{align*}
Then (\ref{eq:lowerBound_2}) follows as above, but with $\bar{p}_{n}(x)=\Delta_{n}\sum_{k=2}^{n}\left(\int_{\R^{d}}p_{t_{k-1}}(x-y)\mu(y)dy\right)$.
By the Plancherel Theorem and $\c Fp_{t_{k-1}}(u)=\c Fp(t_{k-1}^{1/2}u)=e^{-\frac{1}{2}|u|^{2}t_{k-1}}$
this equals
\begin{align*}
\bar{p}_{n}(x) & =(2\pi)^{-d}\Delta_{n}\sum_{k=2}^{n}\left(\int_{\R^{d}}e^{-\frac{1}{2}|u|^{2}t_{k-1}}\c F\mu(u)e^{i\sc ux}du\right).
\end{align*}
We find $\bar{p}_{n}\rightarrow\bar{p}_{\infty}=(2\pi)^{-d}\int_{\R^{d}}(\int_{0}^{T}e^{-\frac{1}{2}|u|^{2}t}dt)\c F\mu(u)e^{i\sc u{\cdot}}du$
$L^{1}(\R^{d})$ and that $\sup_{x\in\R^{d}}\bar{p}_{n}(x)$ is bounded.
As in the proof of Theorem \ref{thm:lowerBound_Hs} conclude by the
dominated convergence theorem that $\E\left[V_{n}\right]\rightarrow\bar{p}_{\infty}(0)\E\left[\int_{\R^{d}}\bar{h}(x)^{2}dx\right]\gtrsim\bar{p}_{\infty}(0)$.
The result follows from applying
\[
\frac{1-e^{-x}}{x}=1-\frac{1}{2}\int_{0}^{x}(1-\frac{r}{x})e^{-r}dr\geq1-\frac{1}{2}\int_{0}^{x}e^{-r}dr\geq\frac{1}{2},\quad x>0,
\]
to $x=-\frac{1}{2}|u|^{2}T$ such that $\int_{0}^{T}e^{-\frac{1}{2}|u|^{2}t}dt=\frac{2}{\left|u\right|^{2}}(1-e^{-\frac{1}{2}|u|^{2}T})\geq\frac{T}{2}$.
\end{proof}

\paragraph*{Acknowledgement}

Support by the DFG Research Training Group 1845 ``Stochastic Analysis
with Applications in Biology, Finance and Physics'' is gratefully
acknowledged.

\bibliographystyle{apalike2}
\bibliography{bibliography}

\end{document}